\documentclass[12,reqno]{amsart}
\usepackage{amsfonts} 
\usepackage{amsmath,amssymb,amsthm, bm}
\usepackage{mathtools}
\usepackage{shuffle}
\usepackage[shortlabels]{enumitem}
\usepackage[alphabetic]{amsrefs}
\usepackage{empheq}
\usepackage{wasysym}
\usepackage{verbatim}
\usepackage{graphicx}
\usepackage[
bookmarks=true,         
bookmarksnumbered=true, 
colorlinks=true, pdfstartview=FitV, linkcolor=blue, citecolor=blue,
urlcolor=blue]{hyperref}
\usepackage{cancel}

\usepackage{graphicx}
\usepackage{amsfonts}
\usepackage{mathrsfs}
\usepackage{fancyhdr}
\usepackage{amsthm}
\usepackage{cases}
\usepackage{amsmath,amssymb,bm}
\usepackage{enumitem}
\usepackage{xcolor}
\usepackage{upgreek}

\usepackage[leftcaption]{sidecap} 
\usepackage{booktabs}
\usepackage{lipsum}
\usepackage{subfigure}
\usepackage{graphicx}
\usepackage{bm}
\usepackage{caption}
\usepackage{mathrsfs}
\usepackage{amsmath}
\usepackage{amsfonts}
\usepackage{caption} 
\usepackage{floatrow}
\newfloatcommand{capbtabbox}{table}[][\FBwidth]
\usepackage{threeparttable,booktabs}
\PassOptionsToPackage{normalem}{ulem}
\usepackage{ulem}
\usepackage{showlabels} 

\providecolor{added}{rgb}{0,0,1}
\providecolor{deleted}{rgb}{1,0,0}

\newtheorem{theorem}{Theorem}[section]
\newtheorem{lemma}[theorem]{Lemma}
\newtheorem{proposition}[theorem]{Proposition}
\newtheorem{remark}[theorem]{Remark}

\theoremstyle{definition}

\allowdisplaybreaks
\numberwithin{equation}{section}

\usepackage[left=1.5in, right=1.5in, top=0.5in, bottom=2in]{geometry}
\begin{document}
\title[Density convergence of spatial average of SWE]{Density convergence of spatial average  of solution to a  one dimensional stochastic wave equation}
 
	 \author[C.Sun]{\sc Chengbo Sun${}^1$ }  
	 \thanks{${}^1$Corresponding author \\ 
	 Y.Hu was supported by   
 NSERC Discovery grant RGPIN
2024-05941  and a centennial  fund of University of Alberta.  
  }  
	 \address{School of Mathematics, Jilin University, Changchun, 130000,   China}
  \email{ scb0927@gmail.com }
\author[Y.Hu]{ \sc   Yaozhong Hu  }

\address{Department of Mathematical and  Statistical
 Sciences, University of Alberta at Edmonton,
Edmonton, Canada, T6G 2G1}
  \email{ yaozhong@ualberta.ca
}

\maketitle

\begin{abstract}
In this paper  we study the spatial averages of the solution of a one-dimensional stochastic wave equation driven by a Gaussian multiplicative noise, which is white in time and has a homogeneous spatial covariance described by the Riesz kernel. 
We establish the rate of convergence for the uniform distance between the density of spatial averages  and the standard normal density. The proof combines  Malliavin calculus with Stein's method for normal approximations. The  key technical challenges lie  in estimating the $L^p$-norm of the second Malliavin derivative and the existence of negative moments of Malliavin covariance matrix.

\smallskip

\end{abstract}
\section{Introduction}
Consider the following one-dimensional stochastic wave equation  
\begin{equation}\label{eqswe}
	\frac{\partial^2u}{\partial t^2}=\frac{\partial^2u}{\partial x^2}+\sigma(u)\dot{W},\quad\quad t\in\mathbb{R}_+, x\in\mathbb{R}
\end{equation}
with initial condition $u(0,x)=1,\frac{\partial}{\partial t}u(0,x)=0 $  (We assume this simple initial condition to focus our study on the stochastic part and other initial conditions can be discussed in the same way), where $\dot{W}$ is a Gaussian noise on some complete probability space $(\Omega, \mathcal{F}, \mathbb{P})$ that is white in time and   has a homogeneous spatial covariance described by the Riesz kernel. That is to say, the covariance of the centered Gaussian noise is given by
\begin{align*}
	\mathbb{E}\left[\dot{W}(t,x)\dot{W}(s,y)\right]=\delta_0(t-s)|x-y|^{-\beta}
\end{align*}
for any $\beta\in(0,1)$, where $\delta_0$ denotes the Dirac delta function at zero. 
This also means that $\dot W$ is white in time and fractional in space with Hurst parameter $H=1-\beta/2$. It is known that there exists a unique continuous   mild solution to \eqref{eqswe}, satisfying 
\begin{equation}
	u(t,x)=1+\int_{0}^{t}\int_{\mathbb{R}}G(t-s,x-y)\sigma(u(s,y))\,W(ds,dy)\, ,
\end{equation}
where $G(t,x)=\frac{1}{2}\mathbf{1}_{\{|x|< t\}}$ denotes the one-dimensional wave kernel. For the diffusion coefficient $\sigma$  we introduce the following assumption:

\noindent   \textbf{(H)} $\sigma:\mathbb{R}\to\mathbb{R}$ is a twice continuously differentiable function so that $\sigma',\sigma''$ are  bounded and $\sigma(x)\ge c>0 $ for all $x\in\mathbb{R}$.

We are interested in the central limit theorems for the density of the spatial averages of solutions. More specifically, 
for an $R>0$ and let  the normalized spatial averages be defined by 
\begin{equation}\label{spatial average}
	F_{R,t}:=\frac{1}{\sigma_{R,t}}\left(\int_{-R}^R[u(t,x)-1]\,dx\right)\,,\ \mathrm{where}\  \sigma_{R,t}^2:=\mathrm{Var}\left(\int_{-R}^Ru(t,x)\,dx\right)\,.
\end{equation}

Our main result in this paper is 

\begin{theorem}\label{main result}  Let Assumption \textbf{(H)}  be satisfied. 
	Let $u=\left\{u(t,x):(t,x)\in\mathbb{R}_+\times\mathbb{R}\right\}$ be the mild solution (namely it satisfies \eqref{eqmild}) to the stochastic wave equation \eqref{eqswe} with initial condition $u(0,x)=1,\frac{\partial}{\partial t}u(0,x)=0 $. 
	Let $F_{R,t}$ be defined as in \eqref{spatial average}. Then, for fixed $t>0$ and $R$ big enough, $F_{R,t}$ has a density $f_{F_{R,t}}$, and there is a constant $c_{t,\beta}$, independent of $R$ so that 
	\begin{equation}
		\sup_{z\in\mathbb{R}}|f_{F_{R,t}}(z)-\phi(z)|\le \begin{cases}
		c_{t,\beta}R^{-\beta/2}\,,&\beta\in(0,\frac{1}{2})\,,\\
		c_{t,\beta}R^{-\beta/2}(\log R)^{1/2}\,,&\beta=\frac{1}{2}\,,\\
		c_{t,\beta}R^{(\beta-1)/2}\,,&\beta\in(\frac{1}{2},1)\,,
	\end{cases}
	\end{equation}
	where $\phi$ is the density of a standard normal distribution on $\mathbb{R}$.
\end{theorem} 
Our assumption that the noise is colored in space (namely, $\beta\in (0, 1)$)  is an important one.  When $\beta =1$ 
(which corresponds  to the  space-time white noise case), the above Theorem \ref{main result}  is no longer valid   because the estimates obtained in Section \ref{estimation of part 1} are no longer applicable to this case. The detailed explanation is provided at in Remark \ref{r.3.3}.

Now let us  give some account on relevant  studies  of the asymptotic 
behavior of the spatial averages of the solution of a stochastic partial differential equation.    Huang, Nualart and Viitasaari \cite{Huang_Nualart_CLTforSHE_2020} were the first  to study spatial averages of solution for    stochastic heat equation,
providing a bound  for the total variation distance $d_{TV}(F_{R,t}, Z) $ between $F_{R,t}$ and a standard normal random variable $Z$. Their approach is the   Malliavin-Stein method. Later on,  several works  have extended this result  to more general  settings of stochastic heat equations (\cite{Balan_Yuan_SHEtimeindependent_2023, Huang_Nualart_SHEcolored_2020,Nualart_Song_Zheng_PAMrough_2020,Nualart_Xia_Zheng_PAMcolored_2022,Nualart_Zheng_SHEtimefractional_2020}). 

The above result has been also extended
to stochastic  wave equations of the form \eqref{eqswe} (\cite{Balan_SPA_2025,Balan_Nualart_Quer-Sardanyons_2022,Balan_Yuan_HAMtimeindependent_2022,Balan_Yuan_HAMtimeindependentrough_2023,Nualart_Zheng_2dSWE_2021, Nualart_Zheng_SWEfractional_2020, 
	Nualart_Zheng_SWEd12_2022}).  The spatial variable  of the equation can be multi-dimensional and the noise structure can also be general. All of above works are about the distance between the distribution of $F_{R,t}$  and a normal one.  There has been no mention of the  density of the law of   $F_{R,t}$. 

On the other hand when the equation is again a  stochastic heat equation driven by   a Gaussian noise which is  white in time and colored in space,  in a series of two papers \cite{Kuzgun_Nualart_densitySHE_2022,
	Kuzgun_Nualart_densitySHEcolored_2024},   
Kuzgun and Nualart  obtain a  rate of convergence in the uniform distance of the density  associated with the normalized average  $F_{R,t}$ to a  normal density. The main  technique  they used is  the Malliavin–Stein method   first introduced in \cite{Hu_Lu_Nualart_density_2014}.  Our main result (Theorem \ref{main result})  is parallel to   the above work. As is well-known in the application of Malliavin calculus to study the density problems,  one  main difficulty is the existence of negative moment of the Malliavin covariance matrix.  For  stochastic heat equations, the solution of the equation 
is always positive and    has itself negative moments.  This characteristic is very helpful  in proving the existence of negative moments in the case of stochastic heat equations   (\cite{Kuzgun_Nualart_densitySHE_2022,Kuzgun_Nualart_densitySHEcolored_2024}).

However, 
in contrast to stochastic heat equations,  for stochastic wave equation \eqref{eqswe} the solution   is no longer  necessarily positive. This makes our task much more challenging. 
Our approach, inspired by the work   \cite{Nualart_density_2007} of  Nualart and Quer-Sardanyons, consists of two key points. The first one is that   we replace $v=-DL^{-1}F$ by $v=DF$
so that the Malliavin covariance matrix $\|DF\|_{\mathcal{H}}^2$ remains non-negative. Second one is that  we assume  the diffusion coefficient  to satisfy  $\sigma(u)\ge c>0$   to  guarantee  the existence of the negative moments of $\|DF\|_{\mathcal{H}}^2$.  A detailed proof of the existence of the negative moments of $\|DF\|_{\mathcal{H}}^2$ is provided in Section \ref{Negative moments}.

It is interesting to point out that our approach  may  also be able to be extended to  higher dimensions. The main difficulty to overcome is to obtain a desirable  bound  for the second order Malliavin derivative  such as the one 
in Proposition  \ref{second derive},  which is more delicate in the higher dimensional case due to the singularity of the 
corresponding  Green's function.   

The remainder of the paper is organized as follows. In Section \ref{preliminaries}, we present several preliminary results on Malliavin calculus, the mild formulation of equation \eqref{eqswe}, and the properties of its Malliavin derivatives. We also collect some auxiliary results needed throughout the paper. 
Due to the specific structure of the wave equation, our arguments deviate   from those used in \cite{Kuzgun_Nualart_densitySHE_2022}.  This is reflected in Section    \ref{proof}, where it  contains the moment estimates of the second Malliavin derivative of the mild solution, the existence of negative moments,  and the other required bounds to complete the proof of Theorem \ref{main result}. Finally, the Appendix includes a technical lemma used repeatedly throughout the paper. In our proof, we do not distinguish between different positive constants $C$.

\section{Preliminaries}	\label{preliminaries}
In this section, we will recall some preliminary notions and results.
 Let $\mathcal{H}_0$ denote the Hilbert space obtained by completing $C_c^\infty(\mathbb{R})$ 
 (the set of all smooth functions from $\mathbb{R}$ to $\mathbb{R}$ with compact support) with respect to the inner product
	\begin{equation*}
		\langle f,g\rangle_{0}
		=
		\int_{\mathbb{R}^{2}} f(x)g(y)|x-y|^{-\beta}\,dx\,dy\,. 
	\end{equation*}
	We use $\|\cdot\|_0$ to denote the norm induced by this inner product, namely
	 $\|f\|_0^2=\langle f,f\rangle_0.$\ \ 
	This space corresponds to the spatial covariance structure of the time-independent noise.

The $L^p(\Omega)$-norm of a real random variable $X$
is denoted by 	  $\|X\|_p$.
\subsection{Burkholder-Davis-Gundy inequality}

We briefly recall the It\^o-Walsh integral and the Burkholder-Davis-Gundy inequality, partly following \cite[Proposition 2.1]{Kuzgun_Nualart_densitySHEcolored_2024}. Let $\mathcal{H}$ be the Hilbert space defined by the inner product
\begin{equation*}
	\begin{split}
		\langle\phi,\psi\rangle_{\mathcal{H}}
		=&\int_{0}^{\infty} \langle \phi(s,\cdot)\,, \psi(s,\cdot)\rangle_0 ds\\
		=&\int_{0}^{\infty}\int_{\mathbb{R}^{2}}\phi(s,x)\psi(s,y)|x-y|^{-\beta}\,dx\,dy\,ds\,.
	\end{split} 
\end{equation*} 
Let $W$ be a Gaussian noise represented by a centered Gaussian family
    	$\left\{ W(\varphi): \varphi \in C_c^\infty(\mathbb{R}_+ \times \mathbb{R}) \right\}$. Assume that its covariance is given by
	\begin{equation*}
		\mathbb{E}\big[W(\varphi)W(\psi)\big]
		= \langle \varphi,\psi\rangle_{\mathcal{H}}\, .
	\end{equation*}	
	Let $(\mathcal{F}_t)_{t\geq 0}$ be the filtration defined by
	$\mathcal{F}_t = \mathcal{F}_t^0 \vee \mathcal{N}$,
	where $\mathcal{F}_t^0$ denotes the $\sigma$-field generated by all random variables
	$W(\varphi)$ such that $\varphi$ is supported in $[0,t]\times \mathbb{R}$, and
	$\mathcal{N}$ is the collection of all $\mathbb{P}$-null sets.
	A random field
   $ X=\{X(t,x):(t,x)\in \mathbb{R}_+\times \mathbb{R}\}$
	is said to be adapted if, for every $(t,x)\in \mathbb{R}_+\times \mathbb{R}$,
	the random variable $X(t,x)$ is $\mathcal{F}_t$-measurable. For any adapted and jointly measurable random field $X$ satisfying
	\begin{equation}
		\mathbb{E}\big[\|X\|_{\mathcal{H}}^2\big]
		=
		\int_0^\infty
		\int_{\mathbb{R}^{2}}
		\mathbb{E}\big[X(s,x)X(s,y)\big]
		|x-y|^{-\beta}\,dx\,dy\,ds
		<\infty\,,
		\label{eq:integrability_condition}
	\end{equation}
	the stochastic integral
	\[
	\int_{\mathbb{R}_+\times \mathbb{R}}
	X(\tau,\xi)\,W(d\tau,d\xi)
	\]
	is well-defined in the It\^o-Walsh sense. Moreover, the following It\^o-Walsh isometry holds:
	\begin{equation}
		\mathbb{E}\left[
		\left|
		\int_{\mathbb{R}_+\times \mathbb{R}}
		X(s,y)\,W(ds,dy)
		\right|^2
		\right]
		=
		\int_0^\infty
		\int_{\mathbb{R}^{2}}
		\mathbb{E}\big[X(s,x)X(s,y)\big]
		|x-y|^{-\beta}\,dx\,dy\,ds \,.
		\label{eq:ito_walsh_isometry}
	\end{equation}
	The following Burkholder--Davis--Gundy type inequality will be used frequently in this work. 
	\begin{proposition}
		Let $p\geq 2$. There exists a constant $C_p>0$ such that, for every adapted and jointly measurable random field $X$ satisfying \eqref{eq:integrability_condition}, and for every $t\geq 0$, one has
		\begin{equation}
			\mathbb{E}\left[
			\left|
			\int_{[0,t]\times\mathbb{R}}
			X(s,y)\,W(ds,dy)
			\right|^p
			\right]
			\leq
			C_p\,
			\mathbb{E}\left[
			\left(
			\int_0^t
			\int_{\mathbb{R}^{2}}
			X(s,y)X(s,x)
			|x-y|^{-\beta}\,dx\,dy\,ds
			\right)^{p/2}
			\right]\,.
			\label{eq:bdg_inequality}
		\end{equation}
\end{proposition}

\subsection{Basic Malliavin calculus} 
We now recall some basic facts on Malliavin calculus associated with the noise process $W$. For further details, we refer to \cite{Hu_BookGaussian_2017} and \cite{Nualart_Malliavin_2006}, which provides a comprehensive account of Malliavin calculus for Gaussian processes. 

We denote by $C_p^{\infty}(\mathbb{R}^n)$ the set of all infinitely   differentiable functions $f:\mathbb{R}^n\to\mathbb{R}$ such that $f$ and all of its partial derivatives have polynomial growth. Let $\mathcal{S}$ denote the class of smooth random variables such that a random variable $F\in\mathcal{S}$ has the form 
\begin{equation}
	F=f(W(h_1),\cdots,W(h_n))\,,
\end{equation}
where $f$ belongs to $ C_p^{\infty}(\mathbb{R}^n)$, $h_1,\cdots,h_n$ are in $\mathcal{H}$, and $n\ge1$. Then the Malliavin derivative of a smooth random variable $F$ is the $\mathcal{H}$-valued random variable given by
\begin{equation}
	DF=\sum_{i=1}^{n}\partial_if(W(h_1),\cdots,W(h_n))\,h_i\,.
\end{equation}
The derivative operator $D$ is closable from $L^p(\Omega)$ into $L^p(\Omega;\mathcal{H})$ for any $p\ge1$ and we let $\mathbb{D}^{1,p}$ be the completion of $\mathcal{S}$ with respect to the norm
\begin{equation*}
	\|F\|_{1,p}=\left(\mathbb{E}[|F|^p]+\mathbb{E}[\|DF\|_{\mathcal{H}}^p]\right)^{1/p}\,.
\end{equation*}
We denote by $\delta$ the adjoint of $D$ given by the duality formula
\begin{equation}\label{eq delta1}
	\mathbb{E}(\delta(u)F)=\mathbb{E}(\langle DF,u\rangle_{\mathcal{H}})
\end{equation}
for any $F\in\mathbb{D}^{1,2}$ and $u\in\mathrm{Dom}\ \delta\subset L^2(\Omega;\mathcal{H}) $, the domain of $\delta$. The operator $\delta$ is also known as the Skorohod integral, as it extends the It\^o integral in the case of Brownian motion.
The following proposition describes a key identity that allows factoring out a scalar random variable inside  a divergence expression.
\begin{proposition}
	Let $F\in\mathbb{D}^{1,2}$ and $u\in\mathrm{Dom}\ \delta$ such that $Fu\in L^2(\Omega;\mathcal{H})$. Then $Fu\in\mathrm{Dom}\ \delta$ and 
	\begin{equation}\label{eq3}
		\delta(Fu)=F\delta(u)-\langle DF,u\rangle_{\mathcal{H}}\,.
	\end{equation}	
\end{proposition}
The following proposition   gives the relationship between the derivative operator $D$ and the divergence operator $\delta$, which will be used in the proof of Theorem \ref{main result}. For additional background, we  refer to \cite[Proposition 6.17]{Hu_BookGaussian_2017} and \cite[Proposition 3.1]{Nualart_Zakai_1986}.  \begin{proposition}
	Let $F\in\mathbb{D}^{1,2}$.  Then
	\begin{align}\label{eq delta2}
		\mathbb{E}[\delta(F)]^2=&\mathbb{E}\left[\|F\|_{\mathcal{H}}^2+\langle D_{\cdot}F_{\star},D_{\star}F_{\cdot}\rangle_{\mathcal{H}^{\otimes 2}}\right]\nonumber\\
		\le&\mathbb{E}\left[\|F\|_{\mathcal{H}}^2\right]+\mathbb{E}\left[\|D_{\cdot}F_{\star}\|_{\mathcal{H}^{\otimes 2}}\right]\,.
	\end{align}
\end{proposition}

\begin{proposition}\label{density}
	Let $F\in \mathbb{D}^{1,1}$ and $v\in L^1(\Omega;\mathcal{H})$ be such that $D_vF\neq0$ a.s. Assume that $v/D_vF\in \mathrm{Dom}\delta$. Then the law of F is absolutely  continuous with respect to the Lebesgue   measure and the density is continuous and bounded  and is  given by
	\begin{equation}
		f_F(x)=\mathbb{E}\left[\mathbf{1}_{\{F>x\}}\delta\left(\frac{v}{D_vF}\right)\right]\,.
	\end{equation}
\end{proposition}

\subsection{Mild solution of stochastic wave equation} \label{Mild solution of SWE}
A random field $(u(t,x), t\in \mathbb{R}_+, x\in \mathbb{R})$  is called a  mild solution to \eqref{eqswe} if it satisfies the following  stochastic integral equation
\begin{equation}\label{eqmild}
	u(t,x)=1+\int_{0}^{t}\int_{\mathbb{R}}G(t-s,x-y)\sigma(u(s,y))\,W(ds,dy)\,,
\end{equation}	 
provided that the above stochastic integral exists, where
\begin{equation}\label{G(t,x)}
    G(t,x)=\frac{1}{2}\mathbf{1}_{\{|x|< t\}}
\end{equation}
is the fundamental solution to the one-dimensional wave equation, and the stochastic integral is understood in the It\^{o}-Walsh sense. Under this noise structure and the Lipschitz condition on the coefficient $\sigma$, it is well-known (see also Proposition \ref{mildsolu}) that equation \eqref{eqswe} admits a unique mild solution
(see e.g.  \cite[Theorem 13]{Dalang}).
More specifically, we have 
\begin{proposition}\label{mildsolu}
	If $\sigma$ is Lipschitz continuous, then \eqref{eqswe} has a unique solution $u(t,x)$, and $u(t,x)$ is  $L^2$-continuous and for any $T>0 $ and $p\ge1$,
	\begin{equation}\label{bound of u}
		\sup_{0\le t\le T}\sup_{x\in\mathbb{R}}\mathbb{E}(|u(t,x)|^p)<\infty\,.
	\end{equation}
\end{proposition}
The following estimate for the moments of the Malliavin derivative of the solution appears in \cite[Lemma 2.2]{Nualart_Zheng_SWEfractional_2020}.
\begin{proposition}\label{p.2.5} 
	Let $u$ be the solution to the stochastic wave equation \eqref{eqswe}.	For any $p\ge2,(t,x)\in\mathbb{R}_+\times\mathbb{R}$, the solution $u(t,x)$ belongs to $\mathbb{D}^{1,p}$. Moreover, for $t> s$, the Malliavin derivative $D_{s,y}u(t,x)$ satisfies
	\begin{align}\label{Malliavin derivate of u}
		D_{s,y}u(t,x)=
		G(t-s,x-y)\sigma(u(s,y))+\int_{s}^{t}\int_{\mathbb{R}}G(t-\tau,x-\xi)\sigma'(u(\tau,\xi))D_{s,y}u(\tau,\xi)\,W(d\tau,d\xi)\,.
	\end{align}
	For any $p\in[2,+\infty),  T >0$ there is a constant   $C=C_{T,p, \beta}$ that depends only on $T,p,\beta $ and the Lipschitz constant of  $\sigma$ so that 
	\begin{equation}\label{ub1}
		\|D_{s,y}u(t,x)\|_p\le CG(t-s,x-y)
	\end{equation}
	for all $0<s<t\le T$ and $x, y\in\mathbb{R}$. 
\end{proposition}
We now present two  auxiliary results from   \cite[Lemma 2.1 and Proposition 3.3]{Nualart_Zheng_SWEfractional_2020} that will be used in the subsequent section.
\begin{lemma}\label{int1}
	Let $G(t,x) $  be defined by \eqref{G(t,x)}. For any $\beta\in(0,1), s,t\ge0$ and $x,\xi\in\mathbb{R}$, we have
	\begin{align}\label{eq2}
		&c_{\beta}\int_{\mathbb{R}^2}G(t,x-y)G(s,\xi-z)|y-z|^{-\beta}\,dy\,dz\nonumber\\
		=&|x-\xi-t-s|^{2-\beta}+|x-\xi+t+s|^{2-\beta}-|x-\xi+t-s|^{2-\beta}-|x-\xi-t+s|^{2-\beta}\,,
	\end{align}
	where $c_{\beta}$ is a constant that depends on $\beta$.
\end{lemma}
\begin{proposition}
	For $\beta\in(0,1)$, denote $\eta(s)=\mathbb{E}[\sigma(u(s,x))]$, which does not depend on $x$ as a consequence of stationary. Then
	\begin{equation}\label{asymptotic of sigma}
		\lim_{R\to\infty}\frac{\sigma_{R,t}^2}{R^{2-\beta}}=2^{2-\beta}\int_{0}^{t}(t-s)^2\eta^2(s)\,ds\,.
	\end{equation}
\end{proposition}
We next present two estimates for integrals involving the one-dimensional wave kernel, which will be used later.
\begin{lemma}
    Let $G(t,x) $  be defined by \eqref{G(t,x)} . For any $\beta\in(0,1), 0<r<s<\tau<t\le T$, define
     \begin{align}\label{LemmaMM}
       M_{t,\beta}
       := &\int_{\mathbb{R}^{2}}G(t-\tau,x-\xi)G(t-\tau,x-\tilde{\xi})G(\tau-r,\xi-z)G(\tau-r,\tilde{\xi}-z)\nonumber\\
		&\quad\times G(\tau-s,\xi-y)G(\tau-s,\tilde{\xi}-y)|\xi-\tilde{\xi}|^{-\beta} \,d\xi \,d\tilde{\xi}\,.
    \end{align}
    We have
    \begin{equation}\label{MM}
        M_{t,\beta}\le c_{T,\beta}G^2(t-s,x-y)G^2(t-r,x-z)\,,
    \end{equation}
    where $c_{T,\beta}$ is a constant that depends on $T$ and $\beta$.
    \begin{proof}
        Since $G(t,x)=\frac{1}{2}\mathbf{1}_{\{|x|< t\}}$ and $0<r<s<\tau<t\le T$, we have
    \begin{align*}
     M_{t,\beta}= & \frac{1}{2^6}\int_{\mathbb{R}^{2}} \mathbf{1}_{\{|x-\xi|< t-\tau\}}\mathbf{1}_{\{|x-\tilde{\xi}|< t-\tau\}}\mathbf{1}_{\{|\xi-y|< \tau-s\}}\mathbf{1}_{\{|\tilde{\xi}-y|< \tau-s\}} \mathbf{1}_{\{|\xi-z|< \tau-r\}}\mathbf{1}_{\{|\tilde{\xi}-z|< \tau-r\}}|\xi-\tilde{\xi}|^{-\beta} \,d\xi \,d\tilde{\xi}\\
     \le&\frac{1}{2^6}\mathbf{1}_{\{|x-y|< t-s\}}\mathbf{1}_{\{|x-z|< t-r\}}\int_{[-t+r,t-r]^{2}}|\xi-\tilde{\xi}|^{-\beta} \,d\xi\, d\tilde{\xi}\\
     =&\frac{2^{1-\beta}}{(1-\beta)(2-\beta)}(t-r)^{2-\beta}G^2(t-s,x-y)G^2(t-r,x-z)\\
     \le&c_{T,\beta}G^2(t-s,x-y)G^2(t-r,x-z)\,,
    \end{align*}
    which is \eqref{MM}. 
    \end{proof}
\end{lemma}
\begin{lemma}
    Let $G(t,x) $  be defined by \eqref{G(t,x)} . For any $\beta\in(0,1), 0<s<\tau<t\le T$, define
     \begin{align}\label{LemmaNN}
       N_{t,\beta}
       := &\int_{\mathbb{R}^{2}}G(t-\tau,x-\xi)G(t-\tau,x-\tilde{\xi})G(\tau-s,\xi-y)G(\tau-s,\tilde{\xi}-\tilde{y})|\xi-\tilde{\xi}|^{-\beta}\, d\xi \,d\tilde{\xi}\,.
    \end{align}
    We have
    \begin{equation}\label{NN}
        N_{t,\beta}\le c_{T,\beta}G(t-s,x-y)G(t-s,x-\tilde{y})\,,
    \end{equation}
    where $c_{T,\beta}$ is the same constant as in \eqref{MM}, depending only on $T$ and $\beta$.
    \begin{proof}
        Since $G(t,x)=\frac{1}{2}\mathbf{1}_{\{|x|< t\}}$ and $0<s<\tau<t\le T$, we have
    \begin{align*}
     N_{t,\beta}= & \frac{1}{2^4}\int_{\mathbb{R}^{2}} \mathbf{1}_{\{|x-\xi|< t-\tau\}}\mathbf{1}_{\{|x-\tilde{\xi}|< t-\tau\}}\mathbf{1}_{\{|\xi-y|< \tau-s\}}\mathbf{1}_{\{|\tilde{\xi}-\tilde{y}|< \tau-s\}} |\xi-\tilde{\xi}|^{-\beta} \,d\xi\, d\tilde{\xi}\\
     \le&\frac{1}{2^4}\mathbf{1}_{\{|x-y|< t-s\}}\mathbf{1}_{\{|x-\tilde{y}|< t-s\}}\int_{[-t+s,t-s]^{2}}|\xi-\tilde{\xi}|^{-\beta}\, d\xi\, d\tilde{\xi}\\
     =&\frac{2^{1-\beta}}{(1-\beta)(2-\beta)}(t-s)^{2-\beta}G(t-s,x-y)G(t-s,x-\tilde{y})\\
     \le&c_{T,\beta}G(t-s,x-y)G(t-s,x-\tilde{y})\,, 
    \end{align*}
    which is \eqref{NN}. 
    \end{proof}
\end{lemma}

\section{Proof  of Theorem \ref{main result}}\label{proof}

\subsection{A different  integration by parts formula for the density}\ 
We shall use an  integration by parts formula different than that   in
\cite[Theorem 3.2]{Kuzgun_Nualart_densitySHE_2022}.  This is necessary since
to prove the existence of the negative moment of the Malliavin covariance matrix Kuzgun and Nualart use the positivity of the solution of the stochastic heat equation which does  not hold true for the solution of a stochastic wave equation. 
The density formula we shall use is  as follows. 

\begin{theorem}\label{density func}
	Assume that $w\in\mathbb{D}^{1,6}(\Omega;\mathcal{H})$, $F=\delta(w)\in\mathbb{D}^{2,6}$ with $\mathbb{E}[F]=0,\mathbb{E}[F^2]=1$, $v=DF$ and $(D_vF)^{-1}\in L^4(\Omega)$. Then $v/D_vF\in \mathrm{Dom} (\delta)$ and $F$ admits a density $f_F(x)$ and the following inequality holds true
	\begin{align}\label{ub of density}
		\sup_{x\in\mathbb{R}}|f_F(x)-\phi(x)|\le&\left\|\|DF\|_{\mathcal{H}}^{-2}\right\|_4\bigg( \left\|\delta(DF)-F\mathbb{E}\left[\|DF\|_{\mathcal{H}}^2\right]\right\|_2+\|F\|_4\sqrt{\mathrm{Var}\left(\|DF\|_{\mathcal{H}}^2\right)}\nonumber\\
		&+2\left\|\|D^2F\|_{\mathcal{H}\otimes\mathcal{H}}\right\|_2\bigg)+2\sqrt{ \mathrm{Var}[\langle DF,w\rangle_{\mathcal{H}}]}\,,
	\end{align}	
	where $\phi(x)$ is the density of the law $N(0,1)$.
\end{theorem}

\begin{proof}
	From Proposition \ref{density},  it follows 
	\begin{equation*}
		f_F(x)=\mathbb{E}\left[\mathbf{1}_{\{F>x\}}\delta\left(\frac{v}{D_vF}\right)\right]\,.
	\end{equation*}
	This implies
	\begin{equation}\label{eq7}
		\sup_{x\in\mathbb{R}}|f_F(x)-\phi(x)|=\sup_{x\in\mathbb{R}}\left|\mathbb{E}\left[\mathbf{1}_{\{F>x\}}\delta\left(\frac{v}{D_vF}\right)\right]-\mathbb{E}\left[\mathbf{1}_{\{N>x\}}N\right]\right|\,,
	\end{equation}	
	where $N\in N(0,1)$ is standard normal random variable.  By \eqref{eq3}, we have
	\begin{equation}\label{eq8}
		\delta\left(\frac{v}{D_vF}\right)=\frac{\delta(v)}{D_vF}-D_v\left(\frac{1}{D_vF}\right)=\frac{\delta(v)}{D_vF}+\frac{D_v(D_vF)}{(D_vF)^2}\,.
	\end{equation}
	Substituting \eqref{eq8} into \eqref{eq7}  yields
	\begin{align*}
		\Phi_x:=&\left|\mathbb{E}\left[\mathbf{1}_{\{F>x\}}\delta\left(\frac{v}{D_vF}\right)\right]-\mathbb{E}\left[\mathbf{1}_{\{N>x\}}N\right]\right|\\
		=&\left|\mathbb{E}\left[\frac{\mathbf{1}_{\{F>x\}}\delta(v)}{D_vF}\right]+\mathbb{E}\left[\frac{\mathbf{1}_{\{F>x\}}D_v(D_vF)}{(D_vF)^2}\right]-\mathbb{E}\left[\mathbf{1}_{\{N>x\}}N\right]\right|\\
		=&\left|\mathbb{E}\left[\frac{\mathbf{1}_{\{F>x\}}\delta(v)}{D_vF}\right]+\mathbb{E}\left[\frac{\mathbf{1}_{\{F>x\}}D_v(D_vF)}{(D_vF)^2}\right]-\mathbb{E}\left[\mathbf{1}_{\{F>x\}}F\right]+\mathbb{E}\left[\mathbf{1}_{\{F>x\}}F\right]-\mathbb{E}\left[\mathbf{1}_{\{N>x\}}N\right]\right|\\
		\le&\mathbb{E}\left[\left|\frac{\delta(v)-FD_vF}{D_vF}\right|\right]+\mathbb{E}\left[\frac{|D_v(D_vF)|}{(D_vF)^2}\right]+|\mathbb{E}\left[F\mathbf{1}_{\{F>x\}}-N\mathbf{1}_{\{N>x\}}\right]|\\
        =:&E_1+E_2+E_3\,.
	\end{align*}
	For the first term $E_1$, by the H\"{o}lder inequality, we have
	\begin{align*}
		&\mathbb{E}\left[\left|\frac{\delta(v)-FD_vF}{D_vF}\right|\right]\\
		=&\mathbb{E}\left[\left|\frac{\delta(v)-F\mathbb{E}\left[\|DF\|_{\mathcal{H}}^2\right]+F\mathbb{E}\left[\|DF\|_{\mathcal{H}}^2\right]-F\|DF\|_{\mathcal{H}}^2}{\|DF\|_{\mathcal{H}}^2}\right|\right]\\
		\le&\mathbb{E}\left[\frac{\left|\delta(v)-F\mathbb{E}\left[\|DF\|_{\mathcal{H}}^2\right]\right|}{\|DF\|_{\mathcal{H}}^2}\right]+\mathbb{E}\left[\frac{\left|F\mathbb{E}\left[\|DF\|_{\mathcal{H}}^2\right]-F\|DF\|_{\mathcal{H}}^2\right|}{\|DF\|_{\mathcal{H}}^2}\right]\\
		\le&\left\|\|DF\|_{\mathcal{H}}^{-2}\right\|_2\left\|\delta(DF)-F\mathbb{E}\left[\|DF\|_{\mathcal{H}}^2\right]\right\|_2+\left\|\|DF\|_{\mathcal{H}}^{-2}\right\|_4\|F\|_4\sqrt{\mathrm{Var}\left(\|DF\|_{\mathcal{H}}^2\right)}\,.
	\end{align*}
	Using the H\"{o}lder inequality to the second term $E_2$ yields
	\begin{align*}
		\mathbb{E}\left[\frac{|D_v(D_vF)|}{(D_vF)^2}\right]&\le\mathbb{E}\left[\|DF\|_{\mathcal{H}}^{-4} \left(\|D^2F\|_{\mathcal{H}\otimes\mathcal{H}}\|DF\|_{\mathcal{H}}^2+\|D^2F\|_{\mathcal{H}\otimes\mathcal{H}}\|DF\|_{\mathcal{H}}^2\right)\right]\\
		&\le 2\left\|\|DF\|_{\mathcal{H}}^{-2}\right\|_2\left\| \|D^2F\|_{\mathcal{H}\otimes\mathcal{H}}\right\|_2\,.
	\end{align*}
	For the third term $E_3$, we apply Stein's method. Specifically,   Stein's equation 
is $f_h'(z)-zf_h(z)=h(z)-\mathbb{E}[h(N)]$
	(e.g. [15, Chapter 3, (3.2.1)]) for the standard normal distribution $N$, whose solution $f_h$ (when  $h(y)=y\mathbf{1}_{\{y>x\}}$)  is given by
	\begin{equation*}
		f_h(x):=e^{x^2/2}\int_{-\infty}^x\{h(y)-\mathbb{E}[h(N)]\}e^{-y^2/2}\,dy\,.
\end{equation*}  
    We have
	\begin{align*}
		|\mathbb{E}\left[F\mathbf{1}_{\{F>x\}}-N\mathbf{1}_{\{N>x\}}\right]|
        &=\left|\mathbb{E}[h(f)-\mathbb{E}[h(N)]]\right|\\
		&=\left|\mathbb{E}[f'_h(F)-Ff_h(F)]\right|\\
		&=\left|\mathbb{E}[f'_h(F)]-\mathbb{E}[\delta(w)f_h(F)]\right|\\
		&=\left|\mathbb{E}[f'_h(F)]-\mathbb{E}[\langle D(f_h(F)),w\rangle_{\mathcal{H}}]\right|\\
		&=\left|\mathbb{E}[f'_h(F)]-\mathbb{E}[f'_h(F)\langle DF,w\rangle_{\mathcal{H}}]\right|\\
		&\le\|f'_h(F)\|_2\|1-\langle DF,w\rangle_{\mathcal{H}}\|_2\,,
	\end{align*}
	where $F=\delta(w)$ and $\|f'_h(F)\|_2\le2$ (see \cite[Lemma 2.1]{Hu_Lu_Nualart_density_2014}), and the last step follows from the Cauchy-Schwarz inequality. Since $1=\mathbb{E}[F^2]=\mathbb{E}[F\delta(w)]=\mathbb{E}[\langle DF,w\rangle_{\mathcal{H}}]$, we have
	\begin{equation*}
		\|1-\langle DF,w\rangle_{\mathcal{H}}\|_2=\sqrt{\mathbb{E}\left[|\mathbb{E}[\langle DF,w\rangle_{\mathcal{H}}]-\langle DF,w\rangle_{\mathcal{H}}|^2\right]}=\sqrt{ \mathrm{Var}[\langle DF,w\rangle_{\mathcal{H}}]}\,.
	\end{equation*}
	Finally, combining the above arguments  we have
	\begin{align*}
		\sup_{x\in\mathbb{R}}|f_F(x)-\phi(x)|\le&\left\|\|DF\|_{\mathcal{H}}^{-2}\right\|_4\bigg( \left\|\delta(DF)-F\mathbb{E}\left[\|DF\|_{\mathcal{H}}^2\right]\right\|_2+\|F\|_4\sqrt{\mathrm{Var}\left(\|DF\|_{\mathcal{H}}^2\right)}\nonumber\\
		&+2\left\|\|D^2F\|_{\mathcal{H}\otimes\mathcal{H}}\right\|_2\bigg)+2\sqrt{ \mathrm{Var}[\langle DF,w\rangle_{\mathcal{H}}]}\,.
	\end{align*}
	This proves the theorem. 
\end{proof}

\subsection{Moment estimates of the second Malliavin derivative of $u$}\label{Moment estimates of the second derivative}
As indicated in Theorem \ref{density func} to prove  Theorem \ref{main result}, we need to complete two tasks: find satisfactory  moment bounds  for the second Malliavin derivative of the solution, and prove the existence of  negative moments of $ \|DF_{R,t}\|_{\mathcal{H}}^2$.  In this subsection we carry out   the first task. 
\begin{proposition}\label{second derive}
	Let $u$ be the solution to the stochastic wave equation \eqref{eqswe}.	For any $p\ge2,(t,x)\in\mathbb{R}_+\times\mathbb{R}$, the solution $u(t,x)$ belongs to $\mathbb{D}^{2,p}$. Moreover, for almost all $0<r<s<t$ and $y,z \in\mathbb{R}$, the second Malliavin derivative $D_{r,z}D_{s,y}u(t,x)$ satisfies
	\begin{align}\label{Second derivative of u}
		D_{r,z}D_{s,y}u(t,x)=&G(t-s,x-y)\sigma'(u(s,y))D_{r,z}u(s,y)\nonumber\\
		&+\int_{[s,t]\times\mathbb{R}}G(t-\tau,x-\xi)\sigma''(u(\tau,\xi))D_{r,z}u(\tau,\xi)D_{s,y}u(\tau,\xi)\,W(d\tau,d\xi)\nonumber\\
		&+\int_{[s,t]\times\mathbb{R}}G(t-\tau,x-\xi)\sigma'(u(\tau,\xi))D_{r,z}D_{s,y}u(\tau,\xi)\,W(d\tau,d\xi)\,.
	\end{align}
	Moreover, for any $p\in[2,+\infty)$  there is  some constant $C=C'_{t,p, \beta}$ that depends on $t,p$ and $\beta$ so that 	\begin{equation}\label{ub2}
		\|D_{r,z}D_{s,y}u(t,x)\|_p\le CG(t-s,x-y)\bigg(G(s-r,y-z)+G(t-r,x-z)\bigg)
	\end{equation}
	for any $0<r<s<t\le T $ and $x,y,z\in\mathbb{R}$.  \end{proposition}
\begin{proof}
	
	For  $0<r<s<t$ and $x,y,z \in\mathbb{R}$,  define 
	$u_0(t,x)=1$,  and for  $ n=0, 1, 2, \cdots, $
	\begin{equation}
		u_{n+1}=1+\int_{0}^{t}\int_{\mathbb{R}}G(t-s,x-y)\sigma(u_n(s,y))\,W(ds,dy)\,.    
	\end{equation} 
	Then it is easy to verify recursively 
	that the 	second Malliavin derivative of $u_n(t,x)$ exists and we have 
	\begin{align}
		D_{r,z}D_{s,y}u_{n+1}(t,x)=&G(t-s,x-y)\sigma'(u_n(s,y))D_{r,z}u_n(s,y)\nonumber\\
		&+\int_{[s,t]\times\mathbb{R}}G(t-\tau,x-\xi)\sigma''(u_n(\tau,\xi))D_{r,z}u_n(\tau,\xi)D_{s,y}u_n(\tau,\xi)\,W(d\tau,d\xi)\nonumber\\
		&+\int_{[s,t]\times\mathbb{R}}G(t-\tau,x-\xi)\sigma'(u_n(\tau,\xi))D_{r,z}D_{s,y}u_n(\tau,\xi)\,W(d\tau,d\xi)\,.\label{e.3.7} 
	\end{align}
	We claim that 
	\begin{equation}\label{claim}
		\|D_{r,z}D_{s,y}u_n(t,x)\|_p\le CG(t-s,x-y)\bigg(G(s-r,y-z)+G(t-r,x-z)\bigg)\,,
	\end{equation}
	for almost all $(s,y),(r,z)\in[0,t]\times\mathbb{R}$, where the constant $C$ may depend on $n$.   We will use the estimate \eqref{ub1} 
	for all $u_n$ with the constant independent of $n$,   which can be proved easily.  We  make use of the mathematical induction. Since $D_{s,y}u_0(t,x)=0$, $D_{s,y}u_1(t,x)=\sigma(1)G(t-s,x-y)$ and $D_{r,z}D_{s,y}u_1(t,x)=0$, we have
	\begin{align*}
		&D_{r,z}D_{s,y}u_{2}(t,x)\\
		=&G(t-s,x-y)\sigma'(u_1(s,y))D_{r,z}u_1(s,y)\\
		&+\int_{[s,t]\times\mathbb{R}}G(t-\tau,x-\xi)\sigma''(u_1(\tau,\xi))D_{r,z}u_1(\tau,\xi)D_{s,y}u_1(\tau,\xi)\,W(d\tau,d\xi)\\
		=&\sigma'(u_1(s,y))\sigma(1)G(t-s,x-y)G(s-r,y-z)\\
		&+\sigma^2(1)\int_{[s,t]\times\mathbb{R}}G(t-\tau,x-\xi)\sigma''(u_1(\tau,\xi))G(\tau-r,\xi-z)G(\tau-s,\xi-y)\,W(d\tau,d\xi)\,.
	\end{align*}
	By assumption \textbf{(H)} , $\sigma(1)$, $\sigma'$ and $\sigma''$ are bounded by a positive constant $L$.
	Using the Burkholder-Davis-Gundy inequality, we have
	\begin{align*}
		&\|D_{r,z}D_{s,y}u_{2}(t,x)\|_p^2\\
		\le&L^4G^2(t-s,x-y)G^2(s-r,y-z)\\
		&+c_pL^6\int_{s}^{t}\int_{\mathbb{R}^{2}}G(t-\tau,x-\xi)G(t-\tau,x-\tilde{\xi})G(\tau-r,\xi-z)G(\tau-r,\tilde{\xi}-z)\\
		&\quad\times G(\tau-s,\xi-y)G(\tau-s,\tilde{\xi}-y)|\xi-\tilde{\xi}|^{-\beta} \,d\xi\, d\tilde{\xi}\,d\tau\\
		=&L^4G^2(t-s,x-y)G^2(s-r,y-z)+c_pL^6M_{t,\beta}\int_s^t\,d\tau\,,
	\end{align*}
    where $M_{t,\beta}$ is defined in \eqref{LemmaMM}.
    By inequality \eqref{MM}, we have
    \begin{align*}
        \|D_{r,z}D_{s,y}u_{2}(t,x)\|_p^2
		\le &L^4G^2(t-s,x-y)G^2(s-r,y-z)\\
        &+c_pc_{T,\beta}L^6(t-s)G^2(t-s,x-y)G^2(t-r,x-z)\,.
    \end{align*}
	So 
	\begin{align*}
		\|D_{r,z}D_{s,y}u_{2}(t,x)\|_p\le C_2G(t-s,x-y)\bigg(G(s-r,y-z)+G(t-r,x-z)\bigg)\,,
	\end{align*} 
	where $C_2=\max\left\{L^2,(c_pc_{T,\beta}t)^{1/2}L^3\right\}$, and $c_p$ is the constant given by the BDG inequality that only depends on $p$. So the claim \eqref{claim} holds true for $n=0,1,2$. Now suppose the claim \eqref{claim} holds true for $n\ge2$.
	By applying the Burkholder-Davis-Gundy inequality, the Minkowski inequality, the assumption \textbf{(H)} and \eqref{ub1},  we have
	\begin{align*}
		&\|D_{r,z}D_{s,y}u_{n+1}(t,x)\|_p^2\\
		\le& L^2C_{T,p,\beta}^2G^2(t-s,x-y)G^2(s-r,y-z)\\
		&+c_pL^2C_{T,p,\beta}^4\int_{s}^{t}\int_{\mathbb{R}^{2}}G(t-\tau,x-\xi)G(t-\tau,x-\tilde{\xi})G(\tau-r,\xi-z)G(\tau-r,\tilde{\xi}-z)\\
		&\quad\times G(\tau-s,\xi-y)G(\tau-s,\tilde{\xi}-y)|\xi-\tilde{\xi}|^{-\beta} \,d\xi\, d\tilde{\xi}\,d\tau\\
		&+c_pL^2C_n^2\int_{s}^{t}\int_{\mathbb{R}^{2}}G(t-\tau,x-\xi)G(t-\tau,x-\tilde{\xi}) G(\tau-s,\xi-y)\bigg(G(s-r,y-z)+G(\tau-r,\xi-z)\bigg)\\
		&\quad\times G(\tau-s,\tilde{\xi}-y)\bigg(G(s-r,y-z)+G(\tau-r,\tilde{\xi}-z)\bigg)|\xi-\tilde{\xi}|^{-\beta} \,d\xi \,d\tilde{\xi}\,d\tau\,,
	\end{align*}
    where the constant $C_{T,p,\beta}$ is given by \eqref{ub1}. Since $G(t,x)=\frac{1}{2}\mathbf{1}_{\{|x|< t\}}$, we have
    \begin{align*}
        &\bigg(G(s-r,y-z)+G(\tau-r,\xi-z)\bigg)\bigg(G(s-r,y-z)+G(\tau-r,\tilde{\xi}-z)\bigg)\\
        =& \frac{1}{2}G(s-r,y-z)+G(s-r,y-z)\bigg(G(\tau-r,\xi-z)+G(\tau-r,\tilde{\xi}-z)\bigg)+G(\tau-r,\xi-z)G(\tau-r,\tilde{\xi}-z)\\
        \le&3G(s-r,y-z)+G(\tau-r,\xi-z)G(\tau-r,\tilde{\xi}-z)\\
        =&6G^2(s-r,y-z)+G(\tau-r,\xi-z)G(\tau-r,\tilde{\xi}-z)\,.
    \end{align*}
    Therefore,
    \begin{align*}
        \|D_{r,z}D_{s,y}u_{n+1}(t,x)\|_p^2
		\le&L^2C_{T,p,\beta}^2G^2(t-s,x-y)G^2(s-r,y-z)+c_pL^2C_{T,p,\beta}^4M_{t,\beta}\int_s^t\,d\tau\\
        &+6c_pL^2C_n^2G^2(s-r,y-z)N_{t,\beta}\int_s^t\,d\tau+c_pL^2C_n^2M_{t,\beta}\int_s^t\,d\tau\,,
    \end{align*}
     where $M_{t,\beta}$ and $N_{t,\beta}$ are given by \eqref{LemmaMM} and \eqref{LemmaNN}, respectively. Then by \eqref{MM} and \eqref{NN}, we have
     \begin{align*}
         &\|D_{r,z}D_{s,y}u_{n+1}(t,x)\|_p^2\\
		\le&L^2C_{T,p,\beta}^2G^2(t-s,x-y)G^2(s-r,y-z)+tc_pL^2C_{T,p,\beta}^4c_{T,\beta}G^2(t-s,x-y)G^2(t-r,x-z)\\
        &+6tc_pL^2C_n^2c_{T,\beta}G^2(t-s,x-y)G^2(s-r,y-z)+tc_pL^2C_n^2c_{T,\beta}G^2(t-s,x-y)G^2(t-r,x-z)\\
        \le &C_{n+1}G(t-s,x-y)\left(G(s-r,y-z)+G(t-r,x-z)\right)\,,
     \end{align*}
	where
	$C_{n+1}=\max\left\{L\left(C_{T,p,\beta}^2+6tc_pC_n^2c_{T,\beta}\right)^{1/2},L\left(tc_pc_{T,\beta}\left(C_{T,p,\beta}^4+C_n^2\right)\right)^{1/2}\right\}$. 
    This shows \eqref{claim} for each $n$ with $C_n$ depending on   $n$.  Next, we want to remove this dependence on  $n$.    
	Iterating \eqref{e.3.7} once we have    
	\begin{align*}
		&D_{r,z}D_{s,y}u_{n+1}(t,x)\\
		=&G(t-s,x-y)\sigma'(u_n(s,y))D_{r,z}u_n(s,y)\\
		&+\int_{[s,t]\times\mathbb{R}}G(t-t_1,x-x_1)\sigma''(u_n(t_1,x_1))D_{r,z}u_n(t_1,x_1)D_{s,y}u_n(t_1,x_1)\,W(dt_1,dx_1)\\
		&+\int_{[s,t]\times\mathbb{R}}G(t-t_1,x-x_1)\sigma'(u_n(t_1,x_1))G(t_1-s,x_1-y)\\
		&\qquad\times \sigma'(u_{n-1}(s,y))D_{r,z}u_{n-1}(s,y)\,W(dt_1,dx_1)\\
		&+\int_{[s,t]\times\mathbb{R}}G(t-t_1,x-x_1)\sigma'(u_n(t_1,x_1))\\
		&\qquad\times\bigg(\int_{[s,t_1]\times\mathbb{R}}G(t_1-t_2,x_1-x_2)\sigma''(u_{n-1}(t_2,x_2))D_{r,z}u_{n-1}(t_2,x_2)\\
		& \qquad \qquad \times D_{s,y}u_{n-1}(t_2,x_2)\,W(dt_2,dx_2)\bigg)\,W(dt_1,dx_1)\\
		&+\int_{[s,t]\times\mathbb{R}}G(t-t_1,x-x_1)\sigma'(u_n(t_1,x_1))\\
		&\qquad\times\bigg(\int_{[s,t_1]\times\mathbb{R}}G(t_1-t_2,x_1-x_2)\sigma' (u_{n-1}(t_2,x_2))\\
		&\qquad\qquad \times D_{r,z}D_{s,y}u_{n-1}(t_2,x_2)\,W(dt_2,dx_2) \bigg)\,W(dt_1,dx_1)\,. 
	\end{align*}
	We continue to iterate to obtain, for $0<r<s<t_n<\cdots<t_2<t_1<t_0=t$ and $x,x_1,\cdots,x_n,y,z\in\mathbb{R}$, 
	\begin{align*}
		&D_{r,z}D_{s,y}u_{n+1}(t,x)\\
		=&G(t-s,x-y)\sigma'(u_n(s,y))D_{r,z}u_n(s,y)\\
		&+\int_{[s,t]\times\mathbb{R}}G(t-t_1,x-x_1)\sigma''(u_n(t_1,x_1))D_{r,z}u_n(t_1,x_1)D_{s,y}u_n(t_1,x_1)\,W(dt_1,dx_1)\\
		&+\int_{[s,t]\times\mathbb{R}}G(t-t_1,x-x_1)\sigma'(u_n(t_1,x_1))G(t_1-s,x_1-y)\\
		&\qquad\times \sigma'(u_{n-1}(s,y))D_{r,z}u_{n-1}(s,y)\,W(dt_1,dx_1)\\
		&+\int_{[s,t]\times\mathbb{R}}G(t-t_1,x-x_1)\sigma'(u_n(t_1,x_1))\\
		&\qquad\times\bigg(\int_{[s,t_1]\times\mathbb{R}}G(t_1-t_2,x_1-x_2)\sigma''(u_{n-1}(t_2,x_2))D_{r,z}u_{n-1}(t_2,x_2)\\
		& \qquad \qquad \times D_{s,y}u_{n-1}(t_2,x_2)\,W(dt_2,dx_2)\bigg)\,W(dt_1,dx_1)\\
		&+\sum_{k=2}^{n}\int_{s}^{t}\cdots\int_{s}^{t_{k-1}}\int_{\mathbb{R}^{k}}G(t_k-s,x_k-y)\sigma'(u_{n-k}(s,y))D_{r,z}u_{n-k}(s,y)\\
		&\qquad\times\prod_{j=1}^kG(t_{j-1}-t_j,x_{j-1}-x_j)\sigma'(u_{n-j+1}(t_j,x_j))\,W(dt_j,dx_j)\\
		&+\sum_{k=2}^{n}\int_{s}^{t}\cdots\int_{s}^{t_{k}}\int_{\mathbb{R}^{k+1}}D_{r,z}u_{n-k+1}(t_k,x_k)D_{s,y}u_{n-k+1}(t_k,x_k)\\
		&\qquad\times\prod_{j=1}^{k-1}G(t_{j-1}-t_j,x_{j-1}-x_j)\sigma'(u_{n-j+1}(t_j,x_j))G(t_{k-1}-t_k,x_{k-1}-x_k)\\
		&\qquad\qquad\times \sigma''(u_{n-k+1}(t_k,x_k))\,W(dt_j,dx_j)\, W(dt_k,dx_k)\\
		=:&\sum_{k=0}^{n}T_k^{(n)}(t,x)+S_k^{(n)}(t,x)\,,
	\end{align*}
	where  $ T_k^{(n)}(t,x)$ and $S_k^{(n)}(t,x)$ denote the $k$th odd and even items in the sum, respectively:
  	\begin{align}
  		S_k^{(n)}(t,x)=&\int_{s}^{t}\cdots\int_{s}^{t_{k}}\int_{\mathbb{R}^{k+1}}D_{r,z}u_{n-k+1}(t_k,x_k)D_{s,y}u_{n-k+1}(t_k,x_k)
        \nonumber \\
  		&\times\prod_{j=1}^{k-1}G(t_{j-1}-t_j,x_{j-1}-x_j)\sigma'(u_{n-j+1}(t_j,x_j))G(t_{k-1}-t_k,x_{k-1}-x_k) \nonumber \\
  		&\quad\times \sigma''(u_{n-k+1}(t_k,x_k))\,W(dt_j,dx_j)\,W(dt_k,dx_k) \nonumber \\
  		=&\int_s^t\int_{\mathbb{R}}G(t-t_1,x-x_1)\sigma'(u_n(t_1,x_1))S_{k-1}^{(n-1)}(t_1,x_1)\,W(dt_1,dx_1) \label{e.3.33} 
  	\end{align} 
    and 
    \begin{align}
       T_k^{(n)}(t,x)=& \int_{s}^{t}\cdots\int_{s}^{t_{k-1}}\int_{\mathbb{R}^{k}}G(t_k-s,x_k-y)\sigma'(u_{n-k}(s,y))D_{r,z}u_{n-k}(s,y)
       \nonumber \\
		&\qquad\times\prod_{j=1}^kG(t_{j-1}-t_j,x_{j-1}-x_j)\sigma'(u_{n-j+1}(t_j,x_j))\,W(dt_j,dx_j)\,,\label{e.3.34} 
    \end{align}
    where $x_0=x,t_0=t$.
    For example,
	\begin{equation*}
		T_0^{(n)}(t,x)=G(t-s,x-y)\sigma'(u_n(s,y))D_{r,z}u_n(s,y)
	\end{equation*}
	and
	\begin{equation*}
		S_0^{(n)}(t,x)=\int_{[s,t]\times\mathbb{R}}G(t-t_1,x-x_1)\sigma''(u_n(t_1,x_1))D_{r,z}u_n(t_1,x_1)D_{s,y}u_n(t_1,x_1)\,W(dt_1,dx_1)\,.
	\end{equation*}
	Now we want to bound all the term $ T_k^{(n)}(t,x)$ and $S_k^{(n)}(t,x)$. 
	
	\noindent \textbf{Case} $k=0$.   It is clear that
	\begin{equation*}
		\|T_0^{(n)}(t,x)\|_p\le  C_{T,p,\beta}LG(t-s,x-y)G(s-r,y-z)\,.
	\end{equation*}
	As for $S_0^{(n)}(t,x)$, applying the Burkholder-Davis-Gundy inequality, the Minkowski inequality, the assumption \textbf{(H)} and \eqref{ub1}, we have
	\begin{align*}
		\|S_0^{(n)}(t,x)\|_p^2\le &c_p\bigg\| \int_{s}^{t}\int_{\mathbb{R}^{2}}G(t-t_1,x-x_1)G(t-t_1,x-\tilde{x}_1)\sigma''(u_n(t_1,x_1))\sigma''(u_n(t_1,\tilde{x}_1))D_{r,z}u_n(t_1,x_1)\\
		&\qquad\times D_{r,z}u_n(t_1,\tilde{x}_1)D_{s,y}u_n(t_1,x_1)D_{s,y}u_n(t_1,\tilde{x}_1)|x_1-\tilde{x}_1|^{-\beta}\,dx_1\,d\tilde{x}_1\,dt_1\bigg\|_{p/2}\\
		\le &c_pL^2 \int_{s}^{t}\int_{\mathbb{R}^{2}}G(t-t_1,x-x_1)G(t-t_1,x-\tilde{x}_1)\|D_{r,z}u_n(t_1,x_1)\\
		&\qquad\times D_{r,z}u_n(t_1,\tilde{x}_1)D_{s,y}u_n(t_1,x_1)D_{s,y}u_n(t_1,\tilde{x}_1)\|_{p/2}|x_1-\tilde{x}_1|^{-\beta}\,dx_1\,d\tilde{x}_1\,dt_1\,.
	\end{align*}
	Using the H\"older inequality and \eqref{MM}, we obtain
	\begin{align*}
		&\|S_0^{(n)}(t,x)\|_p^2\\
		\le&c_pL^2C_{T,2p,\beta}^4 \int_{s}^{t}\int_{\mathbb{R}^{2}}G(t-t_1,x-x_1)G(t-t_1,x-\tilde{x}_1)G(t_1-r,x_1-z)G(t_1-r,\tilde{x}_1-z)\\
		&\times G(t_1-s,x_1-y)G(t_1-s,\tilde{x}_1-y)|x_1-\tilde{x}_1|^{-\beta}\,dx_1\,d\tilde{x}_1\,dt_1\\
        \le &tc_pL^2C_{T,2p,\beta}^4c_{T,\beta}G^2(t-s,x-y)G^2(t-r,x-z)\,.
	\end{align*}
	
	\noindent \textbf{Case} $k=1$.  First, we deal with $T_1^{(n)}(t,x)$. Applying the Burkholder-Davis-Gundy inequality, the Minkowski inequality, the H\"older inequality, the assumption \textbf{(H)} and \eqref{ub1}, we have 
	\begin{align*}
		\|T_1^{(n)}(t,x)\|_p^2\le
		&c_{2p}\bigg\| \int_{s}^{t}\int_{\mathbb{R}^{2}}G(t-t_1,x-x_1)G(t-t_1,x-\tilde{x}_1)\sigma'(u_n(t_1,x_1))\sigma'(u_n(t_1,\tilde{x}_1))G(t_1-s,x_1-y)\\
		&\qquad\times G(t_1-s,\tilde{x}_1-y) |x_1-\tilde{x}_1|^{-\beta}\,dx_1\,d\tilde{x}_1\,dt_1\bigg\|_{p}\|\sigma'(u_{n-1}(s,y))D_{r,z}u_{n-1}(s,y)\|_{2p}^2\\
		\le& c_{2p}L^2C_{T,2p,\beta}^2G^2(s-r,y-z)\int_{s}^{t}\int_{\mathbb{R}^{2}}G(t-t_1,x-x_1)G(t-t_1,x-\tilde{x}_1)G(t_1-s,x_1-y)\\
		&\qquad\times G(t_1-s,\tilde{x}_1-y) \|\sigma'(u_n(t_1,x_1))\sigma'(u_n(t_1,\tilde{x}_1))\|_{p} |x_1-\tilde{x}_1|^{-\beta}\,dx_1\,d\tilde{x}_1\,dt_1\\
		\le& c_{2p}L^4C_{T,2p,\beta}^2G^2(s-r,y-z)\int_{s}^{t}\int_{\mathbb{R}^{2}}G(t-t_1,x-x_1)G(t-t_1,x-\tilde{x}_1)G(t_1-s,x_1-y)\\
		&\qquad\times  G(t_1-s,\tilde{x}_1-y)|x_1-\tilde{x}_1|^{-\beta}\,dx_1\,d\tilde{x}_1\,dt_1\\
        \le&tc_{2p}L^4C_{T,2p,\beta}^2c_{T,\beta}G^2(t-s,x-y)G^2(s-r,y-z)\,,
	\end{align*}
    where the last step follows from \eqref{NN}.
	Now we deal with $\|S_1^{(n)}(t,x)\|_p^2$.  We can write
	\begin{equation*}
		S_1^{(n)}(t,x)=\int_{s}^{t}\int_{\mathbb{R}}G(t-t_1,x-x_1)\sigma'(u_n(t_1,x_1))S_0^{(n-1)}(t_1,x_1)\,W(dt_1,dx_1)
	\end{equation*}
	with
	\begin{equation*}
		S_0^{(n-1)}(t_1,x_1)=\int_{s}^{t_1}\int_{\mathbb{R}}G(t_1-t_2,x_1-x_2)\sigma''(u_{n-1}(t_2,x_2))D_{r,z}u_{n-1}(t_2,x_2)D_{s,y}u_{n-1}(t_2,x_2)\,W(dt_2,dx_2)\,.
	\end{equation*}
	Thus, using the Burkholder-Davis-Gundy inequality, the Minkowski inequality and the assumption \textbf{(H)}, we have
	\begin{align*}
		\|S_1^{(n)}(t,x)\|_p^2\le 
		& c_p\bigg\|\int_{s}^{t}\int_{\mathbb{R}^{2}}G(t-t_1,x-x_1)G(t-t_1,x-\tilde{x}_1)\sigma'(u_n(t_1,x_1))\sigma'(u_n(t_1,\tilde{x}_1))\\
		&\qquad\times S_0^{(n-1)}(t_1,x_1)S_0^{(n-1)}(t_1,\tilde{x}_1)|x_1-\tilde{x}_1|^{-\beta}\,dx_1\,d\tilde{x}_1\,dt_1\bigg\|_{p/2}\\
		\le&c_pL^2\int_{s}^{t}\int_{\mathbb{R}^{2}}G(t-t_1,x-x_1)G(t-t_1,x-\tilde{x}_1)
		\|S_0^{(n-1)}(t_1,x_1)S_0^{(n-1)}(t_1,\tilde{x}_1)\|_{p/2}|x_1-\tilde{x}_1|^{-\beta}\,dx_1\,d\tilde{x}_1\,dt_1\,.
	\end{align*}
	Then using the H\"older inequality, we have
	\begin{align*}
		\|S_1^{(n)}(t,x)\|_p^2\le 
		& c_pL^2\int_{s}^{t}\int_{\mathbb{R}^{2}}G(t-t_1,x-x_1)G(t-t_1,x-\tilde{x}_1)\|S_0^{(n-1)}(t_1,x_1)\|_{p}\|S_0^{(n-1)}(t_1,\tilde{x}_1)\|_{p}
		|x_1-\tilde{x}_1|^{-\beta}\,dx_1\,d\tilde{x}_1\,dt_1\,.
	\end{align*}
	Since  
	\begin{equation*}
		\|S_0^{(n-1)}(t_1,x_1)\|_p\le \sqrt{t_1c_pc_{T,\beta}}LC_{T,2p,\beta}^2G(t_1-s,x_1-y)G(t_1-r,x_1-z)\,,
	\end{equation*}
	by \eqref{MM}, we have
	\begin{align*}
		\|S_1^{(n)}(t,x)\|_p^2\le & c_p^2L^4c_{T,\beta}C_{T,2p,\beta}^4\int_{s}^{t}\int_{\mathbb{R}^{2}}G(t-t_1,x-x_1)G(t-t_1,x-\tilde{x}_1)G(t_1-s,x_1-y)\\
		&\times G(t_1-s,\tilde{x}_1-y)G(t_1-r,x_1-z)G(t_1-r,\tilde{x}_1-z)
		|x_1-\tilde{x}_1|^{-\beta}\,dx_1\,d\tilde{x}_1\,dt_1\\
	\le&c_p^2L^4c_{T,\beta}^2C_{T,2p,\beta}^4G^2(t-s,x-y)G^2(t-r,x-z)\int_S^tt_1\,dt_1\\
    =&\frac{t^2}{2}c_p^2L^4c_{T,\beta}^2C_{T,2p,\beta}^4G^2(t-s,x-y)G^2(t-r,x-z)\,.
	\end{align*}
	
	\noindent \textbf{Case} $2\le k\le n$. We will pay particular attention to the dependence on $k$ and $n$ of the
	constants.  
    Let $t_0=t, x_0=x$.  Using the Burkholder-Davis-Gundy inequality, the Cauchy-Schwartz inequality, the assumption $\sigma'\le L$ in \textbf{(H)},  \eqref{ub1} and the expression \eqref{e.3.34}, we have
	\begin{align*}
		\|T_k^{(n)}(t,x)\|_p^2= &\bigg\|\int_{s}^{t}\cdots\int_{s}^{t_{k-1}}\int_{\mathbb{R}^{k}}G(t_k-s,x_k-y)\sigma'(u_{n-k}(s,y))D_{r,z}u_{n-k}(s,y)\\
		&\quad\times\prod_{j=1}^kG(t_{j-1}-t_j,x_{j-1}-x_j)\sigma'(u_{n-j+1}(t_j,x_j))\,W(dt_j,dx_j)\bigg\|_p^2\\
        \le&L^{2k+2}C_{T,2p,\beta}^2G^2(s-r,y-z)\bigg\|\int_{s}^{t}\cdots\int_{s}^{t_{k-1}}\int_{\mathbb{R}^{k}}G(t_k-s,x_k-y)\\
		&\quad\times\prod_{j=1}^kG(t_{j-1}-t_j,x_{j-1}-x_j)\,W(dt_j,dx_j)\bigg\|_p^2\\
		\le& c_{2p}^kL^{2k+2}C_{T,2p,\beta}^2G^2(s-r,y-z)\int_{s}^{t}\cdots\int_{s}^{t_{k-1}}\int_{\mathbb{R}^{2k}}G(t_k-s,x_k-y)G(t_k-s,\tilde{x}_k-y)\\
		&\quad\times G(t_{j-1}-t_j,x_{j-1}-x_j)G(t_{j-1}-t_j,\tilde{x}_{j-1}-\tilde{x}_j)|x_j-\tilde{x}_j|^{-\beta}dx_kd\tilde{x}_k\cdots \,dx_1\,d\tilde{x}_1\,dt_k\cdots \,dt_1\,.
	\end{align*}
    Using inequality \eqref{NN}, we have
\begin{align*}
	&\int_{\mathbb{R}^{2}}G(t-t_1,x-x_1)G(t-t_1,x-\tilde{x}_1)G(t_1-t_2,x_1-x_2)G(t_1-t_2,\tilde{x}_1-\tilde{x}_2)|x_1-\tilde{x}_1|^{-\beta}\,dx_1\,d\tilde{x}_1\\
	\le&c_{T,\beta}G(t-t_2,x-x_2)G(t-t_2,x-\tilde{x}_2)\,.
\end{align*}
	Then applying \eqref{NN} successively to integrate with respect to $x_2,\tilde{x}_2,\cdots,x_k,\tilde{x}_k$, we have
	\begin{align*}
		\|T_k^{(n)}(t,x)\|_p^2\le & c_{2p}^kc_{T,\beta}^kL^{2k+2}C_{T,2p,\beta}^2G^2(t-s,x-y)G^2(s-r,y-z)\int_{s}^{t}\cdots\int_{s}^{t_{k-1}}\,dt_k\cdots \,dt_1\\
		\le&  \frac{t^k}{k!}c_{2p}^kc_{T,\beta}^kL^{2k+2}C_{T,2p,\beta}^2G^2(t-s,x-y)G^2(s-r,y-z)\,.
	\end{align*}
    As for $\|S_k^{(n)}(t,x)\|_p^2$,  
we shall use mathematical induction on $k$ to show 
	\begin{align}
  		&\|S_{k}^{(n)}(t,x)\|_p^2  \le 
  	 \frac{t^{k+1}}{(k+1)!}c_p^{k+1}L^{2(k+1)}c_{T,\beta}^{k+1}C_{T,2p,\beta}^4G^2(t-s,x-y)G^2(t-r,x-z)\,.
     \label{e.33} 
  \end{align}  
First,  we have 
  	\begin{align*}
  		&\|S_0^{(n)}(t,x)\|_p^2\le tc_pL^2c_{T,\beta}C_{T,2p,\beta}^4G^2(t-s,x-y)G^2(t-r,x-z)\,,\\
  		&\|S_1^{(n)}(t,x)\|_p^2\le \frac{t^2}{2}c_p^2L^4c_{T,\beta}^2C_{T,2p,\beta}^4G^2(t-s,x-y)G^2(t-r,x-z)\,.
  	\end{align*}
 This means 
 that   \eqref{e.33} holds true for $k=0$ and suppose it holds true for $k-1$.  Then using the Burkholder-Davis-Gundy inequality, the Minkowski inequality, the Cauchy-Schwartz inequality, inequality (22), the assumption $\sigma'\le L$, and the expression 
 \eqref{e.3.33}  we have 
  	\begin{align*}
  		&\|S_{k}^{(n)}(t,x)\|_p^2 \\ \le 
  		& c_pL^2\int_{s}^{t}\int_{\mathbb{R}^{2}}G(t-t_1,x-x_1)G(t-t_1,x-\tilde{x}_1)\|S_{k-1}^{(n-1)}(t_1,x_1)\|_{p}\|S_{k-1}^{(n-1)}(t_1,\tilde{x}_1)\|_{p}\\
  		&\quad\times
  		|x_1-\tilde{x}_1|^{-\beta}\,dx_1\,d\tilde{x}_1\,dt_1\\
  		\le&c_p^{k+1}L^{2(k+1)}c_{T,\beta}^{k}C_{T,2p,\beta}^4\int_{s}^{t}\frac{t_1^k}{k!}\int_{\mathbb{R}^{2}}G(t-t_1,x-x_1)G(t-t_1,x-\tilde{x}_1)G(t_1-s,x_1-y)\\
  		&\quad\times G(t_1-r,x_1-z)G(t_1-s,\tilde{x}_1-y)G(t_1-r,\tilde{x}_1-z)
  		|x_1-\tilde{x}_1|^{-\beta}\,dx_1\,d\tilde{x}_1\,dt_1\\
  		\le&\frac{t^{k+1}}{(k+1)!}c_p^{k+1}L^{2(k+1)}c_{T,\beta}^{k+1}C_{T,2p,\beta}^4G^2(t-s,x-y)G^2(t-r,x-z)\,.
  \end{align*} 
	Finally, for any $n\ge3$, we have
\begin{align*}
		&\|D_{r,z}D_{s,y}u_{n+1}(t,x)\|_p\\
		\le &\sum_{k=0}^{n}\|T_k^{(n)}(t,x)\|_p+\|S_k^{(n)}(t,x)\|_p\\
		\le& G(t-s,x-y)\sum_{k=0}^{n}\left(\frac{C_{T,2p,\beta}L}{\sqrt{k!}}(tc_{2p}c_{T,\beta}L^2)^{\frac{k}{2}}G(s-r,y-z)+\frac{C_{T,2p,\beta}^2}{\sqrt{(k+1)!}}(tc_pc_{T,\beta}L^2)^{\frac{k+1}{2}}G(t-r,x-z)\right)\\
		\le&\max\left\{\sum_{k=0}^{n}\frac{C_{T,2p,\beta}L}{\sqrt{k!}}(tc_{2p}c_{T,\beta}L^2)^{\frac{k}{2}},\sum_{k=0}^{n}\frac{C_{T,2p,\beta}^2}{\sqrt{(k+1)!}}(tc_pc_{T,\beta}L^2)^{\frac{k+1}{2}}\right\}\\
        &\times G(t-s,x-y)\bigg(G(s-r,y-z)+G(t-r,x-z)\bigg)\\
		\le& C'_{t,p,\beta}G(t-s,x-y)\bigg(G(s-r,y-z)+G(t-r,x-z)\bigg)\,,
	\end{align*}
	where $C'_{t,p,\beta}=\max\left\{\sum_{k=0}^{n}\frac{C_{T,2p,\beta}L}{\sqrt{k!}}(tc_{2p}c_{T,\beta}L^2)^{\frac{k}{2}},\sum_{k=0}^{n}\frac{C_{T,2p,\beta}^2}{\sqrt{(k+1)!}}(tc_pc_{T,\beta}L^2)^{\frac{k+1}{2}}\right\} $. The constant $C'_{t,p,\beta} $ depends only on $t,T,p$ and $\beta$, and does not depend on $n$.

	Then we need to show that $D_{r,z}D_{s,y}u(t,x) $ is an element of $L^p(\Omega;L^{2}(\mathbb{R}_+\times\mathbb{R})) $. Applying the Minkowski inequality, we have
	\begin{align*}
		&\sup_{n\in\mathbb{N}}\mathbb{E}\left[\|D^2u_n(t,x)\|^p_{\mathcal{H}\otimes\mathcal{H}}\right] \\
		\le&\sup_{n\in\mathbb{N}}\left( \int_0^t\int_{0}^{s}\int_{\mathbb{R}^{2}}\|D_{r,z}D_{s,y}u_{n}(t,x)\|_p\|D_{r,\tilde{z}}D_{s,\tilde{y}}u_{n}(t,x)\|_p|y-\tilde{y}|^{-\beta}|z-\tilde{z}|^{-\beta}\,dz\,d\tilde{z}\,dy\,d\tilde{y}\,dr\,ds\right)^{\frac{p}{2}}\\
		\le &C \bigg( \int_0^t\int_{0}^{s}\int_{\mathbb{R}^{2}}G(t-s,x-y)\bigg(G(s-r,y-z)+G(t-r,x-z)\bigg)\\
        &\quad\times G(t-s,x-\tilde{y})\bigg(G(s-r,\tilde{y}-\tilde{z})+G(t-r,x-\tilde{z})\bigg)|y-\tilde{y}|^{-\beta}|z-\tilde{z}|^{-\beta}\,dz\,d\tilde{z}\,dy\,d\tilde{y}\,dr\,ds\bigg)^{\frac{p}{2}}\\
		\le&C_{\beta}t^{(3-\beta)p}<\infty
	\end{align*}
	uniformly in $n\ge1$ and uniformly in $x\in\mathbb{R}$. In particularly, $\{D^2u_{n+1}(t,x),n\ge1\}$ is uniformly bounded in $L^p(\Omega;L^{2}(\mathbb{R}_+\times\mathbb{R})) $.  This fact together with the $L^p$-convergence of $u_n(t,x)$ to $u(t,x)$   implies (see  \cite[Section 4.1]{Nualart_Zheng_SWEfractional_2020})  the convergence of $D^2u_{n+1}(t,x)$ to $D^2u(t,x)$ in the weak topology on $L^p(\Omega;L^{2}(\mathbb{R}_+\times\mathbb{R})) $ up to a subsequence. In particular, for any fixed $T\in\mathbb{R}_+ $, we have 
	\begin{equation*}
		\sup_{(t,x)\in[0,T]\times\mathbb{R}} \mathbb{E}\left[\|D^2u_n(t,x)\|^p_{\mathcal{H}\otimes\mathcal{H}}\right]  <\infty\,.
	\end{equation*}
	By the same arguments as in the proof of (5.1) in \cite[Section 4.1]{Nualart_Zheng_SWEfractional_2020}, we deduce the estimate
	\begin{equation*}
		\|D_{r,z}D_{s,y}u(t,x)\|_p\le CG(t-s,x-y)\bigg(G(s-r,y-z)+G(t-r,x-z)\bigg)
	\end{equation*}
	for all $0<r<s<t$, $x,y,z\in\mathbb{R}$ and $p\ge2$.
	This completes the proof of the proposition. 
\end{proof}
\subsection{Negative moments $\left\|\|DF_{R, t} \|_{\mathcal{H}}^{-2}\right\|_p$}\label{Negative moments}  In this subsection,  we deal with the negative moment of the Malliavin covariance matrix.   Recall that the  mild solution  $u(t,x)$ to \eqref{eqmild}  and  the random variable $F_{R,t}$ defined  by
\begin{align*}
	F_{R,t}
	&:=\frac{1}{\sigma_{R,t}}\left(\int_{-R}^R[u(t,x)-1]\,dx\right)\\
	&=\frac{1}{\sigma_{R,t}}\left(\int_{-R}^R\int_{0}^{t}\int_{\mathbb{R}}G(t-s,x-y)\sigma(u(s,y))W(ds,dy)\,dx\right)\,,
\end{align*}
where $\sigma_{R,t}^2 
:= \mathrm{Var}\left(\int_{-R}^Ru(t,x)\,dx\right) $.
Using the stochastic Fubini's theorem, we have
\begin{align*}		F_{R,t}&=\int_{0}^{t}\int_{\mathbb{R}}\frac{1}{\sigma_{R,t}}\left(\int_{-R}^RG(t-s,x-y)\sigma(u(s,y))\,dx\right)\,W(ds,dy)\\
	&=\frac{1}{\sigma_{R,t}}\int_{0}^{t}\int_{\mathbb{R}}\phi_{R,t}(s,y)\sigma(u(s,y))\,W(ds,dy)\,,	
\end{align*}
where
\begin{align*}
	\phi_{R,t}(s,y):=\int_{-R}^RG(t-s,x-y)\,dx=\frac{1}{2}\int_{-R}^{R}\mathbf{1}_{\{|x-y|< t-s\}}\,dx\,.
\end{align*}
The Malliavin derivative of $F_{R,t}$   is 
\begin{equation}\label{Malliavin of F_R,t}
	\begin{split} 
		D_{s,y}F_{R,t}=&\int_{-R}^RD_{s,y}u(t,x)\,dx\\
		=& \frac{1}{\sigma_{R,t}}\phi_{R,t}(s,y)\sigma(u(s,y)) +\frac{1}{\sigma_{R,t}}\int_{s}^{t}\int_{\mathbb{R}}\phi_{R,t}(\tau, \xi )D_{s,y} \sigma(u(\tau, \xi))\,W(d\tau, d\xi)\,. 
	\end{split} 
\end{equation}
Let $\delta\in(0,1)$ be arbitrary. We have
\begin{align*}
	\|D_{s,y}F_{R,t}\|_{\mathcal{H}}^2
	=&\int_{0}^{t}\int_{\mathbb{R}^{2}}D_{s,y}F_{R,t}D_{s,\tilde{y}}F_{R,t}|y-\tilde{y}|^{-\beta} \,dy\,d\tilde{y}\,ds\\
	=&\int_{0}^{t}\|D_{s,\cdot}F_{R,t}\|_{0}^2\,ds 
	\ge \int_{t-\delta}^{t}\|D_{s,\cdot}F_{R,t}\|_{0}^2\,ds\,.
\end{align*}
By \eqref{Malliavin derivate of u}, \eqref{Malliavin of F_R,t} and the inequality $\|a+b\|_0^2\ge \frac{1}{2}\|a\|_0^2-\|b\|_0^2$, we have
\begin{align*}
	\|D_{s,y}F_{R,t}\|_{\mathcal{H}}^2
	\ge&\int_{t-\delta}^{t}\|D_{s,\cdot}F_{R,t}\|_{0}^2\,ds\\
	=&\frac{1}{\sigma_{R,t}^2}\int_{t-\delta}^{t}\left\|\phi_{R,t}(s,\cdot)\sigma(u(s,\cdot))+\int_{s}^{t}\int_{\mathbb{R}}\phi_{R,t}(\tau,\xi)\sigma'(u(\tau,\xi)) D_{s,\cdot}u(\tau,\xi)\,W(d\tau,d\xi)\right\|_{0}^2\,ds\\
	\ge&\frac{1}{2\sigma_{R,t}^2}\int_{t-\delta}^{t}\left\|\phi_{R,t}(s,\cdot)\sigma(u(s,\cdot))\right\|_{0}^2\,ds
	-I(\delta)\,.
\end{align*}
where
\begin{equation*}
	I(\delta)=\frac{1}{\sigma_{R,t}^2}\int_{t-\delta}^{t}\left\|\int_{s}^{t}\int_{\mathbb{R}}\phi_{R,t}(\tau,\xi)\sigma'(u(\tau,\xi))D_{s,\cdot}u(\tau,\xi)\,W(d\tau,d\xi)\right\|_{0}^2\,ds\,.
\end{equation*}
Using the assumption that $\sigma(z)\ge c>0$ for all $z\in\mathbb{R}$, we have
\begin{align*}
	&\frac{1}{\sigma_{R,t}^2}\int_{t-\delta}^{t}\left\|\phi_{R,t}(s,\cdot)\sigma(u(s,\cdot))\right\|_{0}^2\,ds\\
	=&\frac{1}{\sigma_{R,t}^2}\int_{t-\delta}^{t}\int_{\mathbb{R}^{2}}\phi_{R,t}(s,y)\phi_{R,t}(s,\tilde{y})\sigma(u(s,y))\sigma(u(s,\tilde{y}))|y-\tilde{y}|^{-\beta}\,dy\,d\tilde{y}\,ds\\
	\ge&c^2g(\delta)\,,
\end{align*}
where
\begin{equation*}
	g(\delta)=\frac{1}{\sigma_{R,t}^2}\int_{t-\delta}^{t}\int_{\mathbb{R}^{2}}\phi_{R,t}(s,y)\phi_{R,t}(s,\tilde{y})|y-\tilde{y}|^{-\beta}\,dy\,d\tilde{y}\,ds\,.
\end{equation*}
Then we can write 
\begin{equation}\label{eq1}
	\int_{0}^{t}\|D_{s,\cdot}F_{R,t}\|_{0}^2\,ds\ge\frac{c^2}{2}g(\delta)-I(\delta)\,.
\end{equation}
So for any $p\in [1, \infty)$, 
\begin{align}\label{eqI_delta1}
	\mathbb{P}(\|D_{s,y}F_{R,t}\|_{\mathcal{H}}^2<\epsilon)
	\le&\mathbb{P}\left(I(\delta)>\frac{c^2}{2}g(\delta)-\epsilon\right)\nonumber\\
	\le&\left( \frac{c^2}{2}g(\delta)-\epsilon\right)^{-p}\mathbb{E}\left[(I(\delta))^p\right]\,.
\end{align} 
We need to find a good upper bound for  $\mathbb{E}
\left[(I(\delta))^p\right]$ and a good lower bound for 
$g(\delta)$. We do the first task first. 
Using the H\"{o}lder inequality, we get
\begin{align*}
	\mathbb{E}\left[(I(\delta))^p\right]
	=&\mathbb{E}\left[\frac{1}{\sigma_{R,t}^2}\int_{t-\delta}^{t}\left\|\int_{s}^{t}\int_{\mathbb{R}}\phi_{R,t}(\tau,\xi)\sigma'(u(\tau,\xi))D_{s,\cdot}u(\tau,\xi)\,W(d\tau,d\xi)\right\|_{0}^2\,ds\right]^p\\		
	\le&\frac{C\delta^{p-1}}{\sigma_{R,t}^{2p}}\mathbb{E}\left[ \int_{t-\delta}^{t}\left\|\int_{s}^{t}\int_{\mathbb{R}}\phi_{R,t}(\tau,\xi)\sigma'(u(\tau,\xi))D_{s,\cdot}u(\tau,\xi)\,W(d\tau,d\xi)\right\|_{0}^{2p}\,ds\right]\,.
\end{align*}
Then we estimate the above expectation. By the Burkholder-Davis-Gundy inequality and the Minkowski inequality, we have
\begin{align}\label{eq4}
	&\mathbb{E}\left[ \left\|\int_{s}^{t}\int_{\mathbb{R}}\phi_{R,t}(\tau,\xi)\sigma'(u(\tau,\xi))D_{s,\cdot}u(\tau,\xi)\,W(d\tau,d\xi)\right\|_{0}^{2p}\right]\nonumber\\
	=&\mathbb{E}\bigg[\bigg|\int_{s}^{t}\int_{s}^{t}\int_{\mathbb{R}^{2}}\int_{\mathbb{R}^{2}}\phi_{R,t}(\tau_1,\xi_1)\phi_{R,t}(\tau_2,\xi_2)\sigma'(u(\tau_1,\xi_1))\sigma'(u(\tau_2,\xi_2))\nonumber\\
	&\quad\times D_{s,y}u(\tau_1,\xi_1)D_{s,\tilde{y}}u(\tau_2,\xi_2)|y-\tilde{y}|^{-\beta}\,dy\,d\tilde{y}\,W(d\tau_1,d\xi_1)\,W(d\tau_2,d\xi_2) \bigg|^p\bigg]\nonumber\\
	\le&c_p\mathbb{E}\bigg[\bigg|\int_{s}^{t}\int_{s}^{t}\int_{\mathbb{R}^{4}}\int_{\mathbb{R}^{2}}\phi_{R,t}(\tau_1,\xi_1)\phi_{R,t}(\tau_1,\tilde{\xi}_1)\phi_{R,t}(\tau_2,\xi_2)\phi_{R,t}(\tau_2,\tilde{\xi}_2)\sigma'(u(\tau_1,\xi_1))\sigma'(u(\tau_2,\xi_2))\nonumber\\
	&\quad\times\sigma'(u(\tau_1,\tilde{\xi}_1))\sigma'(u(\tau_2,\tilde{\xi}_2))D_{s,y}u(\tau_1,\xi_1)D_{s,y}u(\tau_1,\tilde{\xi}_1)D_{s,\tilde{y}}u(\tau_2,\xi_2)D_{s,\tilde{y}}u(\tau_2,\tilde{\xi}_2)\nonumber\\
	&\quad\quad\times|\xi_1-\tilde{\xi}_1|^{-\beta}|\xi_2-\tilde{\xi}_2|^{-\beta}|y-\tilde{y}|^{-\beta}\,dy\,d\tilde{y}\,d\xi_1\,d\tilde{\xi}_1\,d\xi_2\,d\tilde{\xi}_2\,d\tau_1\,d\tau_2\bigg|^{\frac{p}{2}}\bigg]\nonumber\\
	\le&c_p\bigg(\int_{s}^{t}\int_{s}^{t}\int_{\mathbb{R}^{4}}\int_{\mathbb{R}^{2}}\phi_{R,t}(\tau_1,\xi_1)\phi_{R,t}(\tau_1,\tilde{\xi}_1)\phi_{R,t}(\tau_2,\xi_2)\phi_{R,t}(\tau_2,\tilde{\xi}_2)\left\|X_{\tau_1,\xi_1,\tilde{\xi}_1,\tau_2,\xi_2,\tilde{\xi}_2,y,\tilde{y}}\right\|_{p/2}\nonumber\\
	&\quad\times |\xi_1-\tilde{\xi}_1|^{-\beta}|\xi_2-\tilde{\xi}_2|^{-\beta}|y-\tilde{y}|^{-\beta}\,dy\,d\tilde{y}\,d\xi_1\,d\tilde{\xi}_1\,d\xi_2\,d\tilde{\xi}_2\,d\tau_1\,d\tau_2\bigg)^{\frac{p}{2}}\,,
\end{align}
where
\begin{align*}
	X_{\tau_1,\xi_1,\tilde{\xi}_1,\tau_2,\xi_2,\tilde{\xi}_2,y,\tilde{y}}:=&\sigma'(u(\tau_1,\xi_1))\sigma'(u(\tau_2,\xi_2))\sigma'(u(\tau_1,\tilde{\xi}_1))\sigma'(u(\tau_2,\tilde{\xi}_2))D_{s,y}u(\tau_1,\xi_1)\\
	&\times D_{s,y}u(\tau_1,\tilde{\xi}_1)D_{s,\tilde{y}}u(\tau_2,\xi_2)D_{s,\tilde{y}}u(\tau_2,\tilde{\xi}_2).
\end{align*}
By \eqref{ub1}  we have
\begin{align*}
	&\left\|X_{\tau_1,\xi_1,\tilde{\xi}_1,\tau_2,\xi_2,\tilde{\xi}_2,y,\tilde{y}}\right\|_{p/2}\\
	\le& CG(\tau_1-s,\xi_1-y)G(\tau_1-s,\tilde{\xi}_1-y)G(\tau_2-s,\xi_2-\tilde{y})G(\tau_2-s,\tilde{\xi}_2-\tilde{y})\,.
\end{align*}
Next, we carry out the integration in \eqref{eq4} with respect to $y$ and $\tilde{y}$, we have
\begin{align*}
	&\mathbb{E}\bigg[ \Big\|\int_{s}^{t}\int_{\mathbb{R}}\phi_{R,t}(\tau,\xi)\sigma'(u(\tau,\xi))D_{s,\cdot}u(\tau,\xi)\,W(d\tau,d\xi)\Big\|_{0}^{2p}\bigg]\\
	\le&c_pC\bigg(\int_{s}^{t}\int_{s}^{t}\int_{\mathbb{R}^{4}}\phi_{R,t}(\tau_1,\xi_1)\phi_{R,t}(\tau_1,\tilde{\xi}_1)\phi_{R,t}(\tau_2,\xi_2)\phi_{R,t}(\tau_2,\tilde{\xi}_2) |\xi_1-\tilde{\xi}_1|^{-\beta}\\
	&\qquad \times|\xi_2-\tilde{\xi}_2|^{-\beta}\,d\xi_1\,d\tilde{\xi}_1\,d\xi_2\,d\tilde{\xi}_2\,d\tau_1\,d\tau_2\bigg)^{\frac{p}{2}}\\
	=&c_pC\left(\int_{s}^{t}\int_{\mathbb{R}^{2}}\phi_{R,t}(\tau,\xi)\phi_{R,t}(\tau,\tilde{\xi})|\xi-\tilde{\xi}|^{-\beta} \,d\xi\, d\tilde{\xi}\,d\tau\right)^p\,.
\end{align*}
So for $\mathbb{E}\left[(I(\delta))^p\right]$, we have
\begin{align}\label{eqI_delta2}
	\mathbb{E}\left[(I(\delta))^p\right]\le& \frac{C\delta^{p-1}}{\sigma_{R,t}^{2p}}\int_{t-\delta}^{t}\left(\int_{s}^{t}\int_{\mathbb{R}^{2}}\phi_{R,t}(\tau,\xi)\phi_{R,t}(\tau,\tilde{\xi})|\xi-\tilde{\xi}|^{-\beta}\, d\xi \,d\tilde{\xi}\,d\tau\right)^p\,ds\nonumber\\
	\le&\frac{C\delta^{p-1}}{\sigma_{R,t}^{2p}}\left(\int_{t-\delta}^{t}\int_{\mathbb{R}^{2}}\phi_{R,t}(\tau,\xi)\phi_{R,t}(\tau,\tilde{\xi})|\xi-\tilde{\xi}|^{-\beta} \,d\xi \,d\tilde{\xi}\,d\tau\right)^p\nonumber\\
	=&C\delta^{p-1}\left( g(\delta)\right)^p\,.
\end{align}
Combining \eqref{eq1}, \eqref{eqI_delta1} and \eqref{eqI_delta2} yields
\begin{align}
	\mathbb{P}(\|D_{s,y}F_{R,t}\|_{\mathcal{H}}^2<\epsilon)
	\le&C\delta^{p-1}(g(\delta))^p\left( \frac{c^2}{2}g(\delta)-\epsilon\right)^{-p}\,.\label{e.3.12} 
\end{align}

Now we obtain   the lower bound of $g(\delta)$. Recall $\phi_{R,t}(s,y):=\int_{-R}^RG(t-s,x-y)dx$ and using the Fourier transform, we have
\begin{align*}
	g(\delta)=&\frac{1}{\sigma_{R,t}^2}\int_{t-\delta}^{t}\int_{\mathbb{R}^{2}}\phi_{R,t}(s,y)\phi_{R,t}(s,\tilde{y})|y-\tilde{y}|^{-\beta}\,dy\,d\tilde{y}\,ds\\
	=&\frac{c_{d,\beta}}{\sigma_{R,t}^2}\int_{t-\delta}^{t}\int_{[-R,R]^2}\int_{\mathbb{R}}|\xi|^{\beta-1}e^{-i(x_1-x_2)\cdot\xi}\left|\hat{G}_{t-s}(\xi)\right|^2\,d\xi \,dx_1\,dx_2\,ds\\
	=&\frac{c_{d,\beta}}{\sigma_{R,t}^2}\int_{0}^{\delta}\int_{[-R,R]^2}\int_{\mathbb{R}}|\xi|^{\beta-1}e^{-i(x_1-x_2)\cdot\xi}\left|\hat{G}_{s}(\xi)\right|^2\,d\xi\, dx_1\,dx_2\,ds\\
	=& \frac{c_{d,\beta}}{\sigma_{R,t}^2}\int_{0}^{\delta}\int_{\mathbb{R}}|\xi|^{\beta-1}\left|\hat{G}_{s}(\xi)\right|^2\left(\int_{[-R,R]^2}e^{-i(x_1-x_2)\cdot\xi}\,dx_1\,dx_2\right)\,d\xi\, ds\\
	=& \frac{c_{d,\beta}}{\sigma_{R,t}^2}\int_{0}^{\delta}\int_{\mathbb{R}}|\xi|^{\beta-1}\left|\hat{G}_{s}(\xi)\right|^2\left|\int_{-R}^{R}e^{-ix\xi}\,dx\right|^2\,d\xi \,ds\\
	=& \frac{c_{d,\beta}}{\sigma_{R,t}^2}\int_{0}^{\delta}\int_{\mathbb{R}}|\xi|^{\beta-1}\left|\hat{G}_{s}(\xi)\right|^2\frac{4\sin^2(R\xi)}{\xi^2}\,d\xi\, ds\,.
\end{align*}
Using \cite[Appendix A]{Quer-Sardanyons_Sanz-Solé_3Dcontinuity_2004}, for any $\delta\ge0$, we have 
\begin{equation*}
	\int_{0}^{\delta}\left|\hat{G}_{s}(\xi)\right|^2\,ds\ge (2^4\pi^2)^{-1}(\delta\wedge \delta^3)\frac{1}{1+|\xi|^2}\,.
\end{equation*}
Therefore, for sufficiently small $\delta$, we get
\begin{align*}
	g(\delta)\ge&\frac{c_{d,\beta}\delta^3}{ 2^4\pi^2 \sigma_{R,t}^2}\int_{\mathbb{R}}\frac{|\xi|^{\beta-1}}{1+|\xi|^2}\frac{4\sin^2(R\xi)}{\xi^2}\,d\xi\, .
\end{align*}
By the change of variables $\eta=R\xi$, we have
\begin{align*}
	g(\delta)\ge&\frac{c_{d,\beta}\delta^3}{ 2^4\pi^2 \sigma_{R,t}^2}\int_{\mathbb{R}}\frac{|\xi|^{\beta-1}}{1+|\xi|^2}\frac{4\sin^2(R\xi)}{\xi^2}\,d\xi \\
    =&\frac{4c_{d,\beta}\delta^3R^{2-\beta}}{ 2^4\pi^2 \sigma_{R,t}^2}\int_{\mathbb{R}}\frac{|\eta|^{\beta-3}\sin^2\eta}{1+|\eta|^2/R^2}\,d\eta\,.
\end{align*}
Since $R\ge1$, we have
\begin{align*}
	g(\delta)\ge&\frac{4c_{d,\beta}\delta^3R^{2-\beta}}{ 2^4\pi^2 \sigma_{R,t}^2}\int_{\mathbb{R}}\frac{|\eta|^{\beta-3}\sin^2\eta}{1+|\eta|^2}\,d\eta\\
    \ge&\frac{C_{d,\beta}\delta^3R^{2-\beta}}{\sigma_{R,t}^2} 
    \ge  C_{t,d,\beta}\delta^3\,,
\end{align*}
where the constant $C_{t,d,\beta}$ is independent of  $R$ and the last inequality follows from \eqref{asymptotic of sigma}.
At this point, we choose $\delta=\delta(\epsilon)$ in such a way that $g(\delta)=4\epsilon/c^2$ since $g(\delta)\to 0$ as $\delta\to 0$.  The above inequality  $g(\delta)\ge C\delta^3$ gives  $\delta\le C\epsilon^{1/3}$. Hence, by \eqref{e.3.12} we have 
\begin{equation}
	\mathbb{P}(\|D_{s,y}F_{R,t}\|_{\mathcal{H}}^2<\epsilon)\le C\epsilon^{\frac{p-1}{3}}\,.	
\end{equation}
Finally,  for any $q>1$,   using this estimate and the formula $\mathbb{E}[X^{-q}]=q\int_{0}^{\infty}\epsilon^{-q-1}\mathbb{P}(X<\epsilon)\,d\epsilon$, we have
\begin{align}
	\sup_{R\ge1}\mathbb{E}\left[\|D_{s,y}F_{R,t}\|_{\mathcal{H}}^{-2q} \right]&=	\sup_{R\ge1}q\int_{0}^{\infty}\epsilon^{-q-1}\mathbb{P}(\|D_{s,y}F_{R,t}\|_{\mathcal{H}}^2<\epsilon)\,d\epsilon \nonumber \\
	&\le 1+\sup_{R\ge1}q\int_{0}^{1}\epsilon^{-q-1}\mathbb{P}(\|D_{s,y}F_{R,t}\|_{\mathcal{H}}^2<\epsilon)\,d\epsilon\nonumber \\
	&\le 1+Cq\int_{0}^{1}\epsilon^{-q-1}\epsilon^{\frac{p-1}{3}}\,d\epsilon<\infty\label{e.3.14}
\end{align}
if  $p-1\ge 3q$, which is possible  since we can allow $p$ to be arbitrary.

\subsection{Estimation of $\left\|\|D_{r,z}D_{s,y}F_{R,t}\|_{\mathcal{H}\otimes\mathcal{H}}\right\|_2$}\label{estimation of part 3}
Recall  that  the second Malliavin derivative $D_{r,z}D_{s,y}u(t,x)$ satisfies Proposition \ref{second derive},
\begin{align*}
	D_{r,z}D_{s,y}u(t,x)=&G(t-s,x-y)\sigma'(u(s,y))D_{r,z}u(s,y)\nonumber\\
	&+\int_{[s,t]\times\mathbb{R}}G(t-\tau,x-\xi)\sigma''(u(\tau,\xi))D_{r,z}u(\tau,\xi)D_{s,y}u(\tau,\xi)\,W(d\tau,d\xi)\nonumber\\
	&+\int_{[s,t]\times\mathbb{R}}G(t-\tau,x-\xi)\sigma'(u(\tau,\xi))D_{r,z}D_{s,y}u(\tau,\xi)\,W(d\tau,d\xi)\,.
\end{align*}
Applying the Burkholder-Davis-Gundy inequality, stochastic Fubini's theorem,  triangle inequality, the hypothesis \textbf{(H)}, and the estimates \eqref{ub1} and \eqref{ub2}, for all $(t,x)\in[0,T]\times\mathbb{R}$, we have
\begin{align*}
	&\left\|\|D_{r,z}D_{s,y}F_{R,t}\|_{\mathcal{H}\otimes\mathcal{H}}\right\|_2\\
	\le&\frac{C}{\sigma_{R,t}^2}\int_{0}^{t}\int_{0}^{s}\int_{\mathbb{R}^{4}}\phi_{R,t}(s,y)\phi_{R,t}(s,y')G(s-r,y-z)G(s-r,y'-z')\\
	&\qquad\times|z-z'|^{-\beta}|y-y'|^{-\beta}\,dz\,dz'\,dy\,dy'\,dr\,ds\\
	&+\frac{C}{\sigma_{R,t}^2}\int_{0}^{t}\int_{0}^{s}\int_{s}^{t}\int_{\mathbb{R}^6}\phi_{R,t}(\tau,\xi)\phi_{R,t}(\tau,\xi')G(\tau-r,\xi-z)G(\tau-r,\xi'-z')G(\tau-s,\xi-y)\\
	&\qquad\times G(\tau-s,\xi'-y')|\xi-\xi'|^{-\beta}|z-z'|^{-\beta}|y-y'|^{-\beta}\,d\xi\, d\xi'\,dz\,dz'\,dy\,dy'\,d\tau \,dr\,ds\\
	&+\frac{C}{\sigma_{R,t}^2}\int_{0}^{t}\int_{0}^{s}\int_{s}^{t}\int_{\mathbb{R}^6}\phi_{R,t}(\tau,\xi)\phi_{R,t}(\tau,\xi')G(\tau-s,\xi-y)G(\tau-s,\xi'-y')G(s-r,y-z)\\
	&\qquad\times G(s-r,y'-z')|\xi-\xi'|^{-\beta}|z-z'|^{-\beta}|y-y'|^{-\beta}\,d\xi \,d\xi'\,dz\,dz'\,dy\,dy'\,d\tau \,dr\,ds\\
	=&C\Phi_{R,t}+C\Phi_{R,t}^{(1)}+C\Phi_{R,t}^{(2)}\,,
\end{align*}
where
\begin{align}\label{Phi1}
   \Phi_{R,t}:=& \frac{1}{\sigma_{R,t}^2}\int_{0}^{t}\int_{0}^{s}\int_{\mathbb{R}^{4}}\phi_{R,t}(s,y)\phi_{R,t}(s,y')G(s-r,y-z)G(s-r,y'-z')\nonumber\\
	&\qquad\times|z-z'|^{-\beta}|y-y'|^{-\beta}\,dz\,dz'\,dy\,dy'\,dr\,ds\,,
\end{align}
\begin{align}\label{Phi2}
    \Phi_{R,t}^{(1)}:=&\frac{1}{\sigma_{R,t}^2}\int_{0}^{t}\int_{0}^{s}\int_{s}^{t}\int_{\mathbb{R}^6}\phi_{R,t}(\tau,\xi)\phi_{R,t}(\tau,\xi')G(\tau-r,\xi-z)G(\tau-r,\xi'-z')G(\tau-s,\xi-y)\nonumber\\
	&\qquad\times G(\tau-s,\xi'-y')|\xi-\xi'|^{-\beta}|z-z'|^{-\beta}|y-y'|^{-\beta}\,d\xi \,d\xi'\,dz\,dz'\,dy\,dy'\,d\tau \,dr\,ds\,,
\end{align}
and
\begin{align}\label{Phi3}
   \Phi_{R,t}^{(2)}:=& \frac{1}{\sigma_{R,t}^2}\int_{0}^{t}\int_{0}^{s}\int_{s}^{t}\int_{\mathbb{R}^6}\phi_{R,t}(\tau,\xi)\phi_{R,t}(\tau,\xi')G(\tau-s,\xi-y)G(\tau-s,\xi'-y')G(s-r,y-z)\nonumber\\
	&\qquad\times G(s-r,y'-z')|\xi-\xi'|^{-\beta}|z-z'|^{-\beta}|y-y'|^{-\beta}\,d\xi \,d\xi'\,dz\,dz'\,dy\,dy'\,d\tau \,dr\,ds\,.
\end{align}
Combining the estimates \eqref{Phi}, \eqref{Phi_2} and \eqref{Phi_3} for $\Phi_{R,t}$, $\Phi_{R,t}^{(1)}$ and $\Phi_{R,t}^{(2)}$, respectively, we obtain
\begin{align}\label{second detivative estimate}
	\left\|\|D_{r,z}D_{s,y}F_{R,t}\|_{\mathcal{H}\otimes\mathcal{H}}\right\|_2\le \begin{cases}
		c_{t,\beta}R^{-\beta}\,,&\beta\in(0,\frac{1}{2})\,,\\
		c_{t,\beta}R^{\beta-1}\log R\,,&\beta=\frac{1}{2}\,,\\
		c_{t,\beta}R^{\beta-1}\,,&\beta\in(\frac{1}{2},1)\,.
	\end{cases}
\end{align}

\subsection{Estimation of $\left\|\delta(D_{s,y}F_{R,t})-F_{R,t}\mathbb{E}\left[\|D_{s,y}F_{R,t}\|_{\mathcal{H}}^2\right]\right\|_2$}\label{estimation of part 1}
Denote   $h(R):=\mathbb{E}\left[\|D_{s,y}F_{R,t}\|_{\mathcal{H}}^2\right]$, we have
\begin{align*}
	&\left\|\delta(D_{s,y}F_{R,t})-F_{R,t}\mathbb{E}\left[\|D_{s,y}F_{R,t}\|_{\mathcal{H}}^2\right]\right\|_2^2 \\
	=&\left\|\delta(D_{s,y}F_{R,t})-F_{R,t}h(R) \right\|_2^2 \\   
	=& \mathbb{E}\left[\left[\delta (D_{s,y}F_{R,t})\right]^2 \right]-2h(R)\mathbb{E}\left[\delta(D_{s,y}F_{R,t})F_{R,t}\right]+h(R)^2\mathbb{E}\left[F_{R,t}^2\right]\,.
\end{align*} 
Using the results \eqref{eq delta1}, \eqref{eq delta2} and the fact that $\mathbb{E}\left[F_{R,t}^2\right]=1$, we have
\begin{align*}
	&\left\|\delta(D_{s,y}F_{R,t})-F_{R,t}\mathbb{E}\left[\|D_{s,y}F_{R,t}\|_{\mathcal{H}}^2\right]\right\|_2^2\\
	 \le& \mathbb{E}\left[\|D_{s,y}F_{R,t}\|_{\mathcal{H}}^2\right]+ \mathbb{E}\left[\| D_{r,z}D_{s,y}F_{R,t}\|_{\mathcal{H}^{\otimes 2}}\right] -2h(R)\mathbb{E}\left[\|D_{s,y}F_{R,t}\|_{\mathcal{H}}^2\right]+h(R)^2\\
	 =&h(R)-h(R)^2+\mathbb{E}\left[\| D_{r,z}D_{s,y}F_{R,t}\|_{\mathcal{H}^{\otimes 2}}\right]\,.
\end{align*}
By the Cauchy-Schwarz inequality and the estimate \eqref{second detivative estimate}, we have
\begin{align*}
	\mathbb{E}\left[\| D_{r,z}D_{s,y}F_{R,t}\|_{\mathcal{H}^{\otimes 2}}\right]\le&\left(\mathbb{E}\left[\left|\| D_{r,z}D_{s,y}F_{R,t}\|_{\mathcal{H}^{\otimes 2}}\right|^2\right]\right)^{\frac{1}{2}}\le \begin{cases}
		c_{t,\beta}R^{-\beta}\,,&\beta\in(0,\frac{1}{2})\,,\\
		c_{t,\beta}R^{\beta-1}\log R\,,&\beta=\frac{1}{2}\,,\\
		c_{t,\beta}R^{\beta-1}\,,&\beta\in(\frac{1}{2},1)\,.
	\end{cases}
\end{align*}
As for $\mathbb{E}\left[\|D_{s,y}F_{R,t}\|_{\mathcal{H}}^2\right] $, recall that
\begin{align}
	D_{s,y}F_{R,t}=\frac{	\phi_{R,t}(s,y)\sigma(u(s,y))}{\sigma_{R,t}}
	+\frac{1}{\sigma_{R,t}}\int_{s}^{t}\int_{\mathbb{R}}\phi_{R,t}(\tau,\xi)\sigma'(u(\tau,\xi))D_{s,y}u(\tau,\xi)\,W(d\tau,d\xi)\,,
\end{align}
then we have
\begin{align*}
	\mathbb{E}\left[\|D_{s,y}F_{R,t}\|_{\mathcal{H}}^2\right] =&\mathbb{E}\left[\frac{1}{\sigma_{R,t}^2} \left\| \phi_{R,t}(s,y)\sigma(u(s,y))\right\|_{\mathcal{H}}^2 \right]\\
	& +\mathbb{E}\left[ \frac{1}{\sigma_{R,t}^2}\left\| \int_{s}^{t}\int_{\mathbb{R}}\phi_{R,t}(\tau,\xi)\sigma'(u(\tau,\xi))D_{s,y}u(\tau,\xi)\,W(d\tau,d\xi)\right\| _{\mathcal{H}}^2\right]\\
	=&:J_1+J_2\,.
\end{align*}
Since  
\begin{align*}
	\sigma_{R,t}^2&=\mathrm{Var}\left(\int_{0}^{t}\int_{\mathbb{R}}\phi_{R,t}(s,y)\sigma(u(s,y))\,W(ds,dy)\right)\\
	&=\int_{0}^{t}\int_{\mathbb{R}^{2}}\phi_{R,t}(s,y)\phi_{R,t}(s,\tilde{y})\mathbb{E}\left[\sigma(u(s,y))\sigma(u(s,\tilde{y}))\right]|y-\tilde{y}|^{-\beta}\,dy\,d\tilde{y}\,ds\,.
\end{align*}
We have
\begin{align*}
	J_1=&\mathbb{E}\left[ \frac{1}{\sigma_{R,t}^2}\left\| \phi_{R,t}(s,y)\sigma(u(s,y))\right\|_{\mathcal{H}}^2 \right]\\
	=&\frac{1}{\sigma_{R,t}^2}\int_{0}^{t}\int_{\mathbb{R}^{2}}\phi_{R,t}(s,y)\phi_{R,t}(s,\tilde{y})\mathbb{E}\left[\sigma(u(s,y))\sigma(u(s,\tilde{y}))\right]|y-\tilde{y}|^{-\beta}\,dy\,d\tilde{y}\,ds\\
	=&1\,.                             
\end{align*}
By the Burkholder-Davis-Gundy inequality and the Minkowski inequality, we have
\begin{align*}
	J_2=&\mathbb{E}\left[ \frac{1}{\sigma_{R,t}^2}\left\| \int_{s}^{t}\int_{\mathbb{R}}\phi_{R,t}(\tau,\xi)\sigma'(u(\tau,\xi))D_{s,y}u(\tau,\xi)\,W(d\tau,d\xi)\right\| _{\mathcal{H}}^2\right]\\
	\le&\frac{C}{\sigma_{R,t}^2}\mathbb{E}\bigg[ \int_{0}^{t}\int_{0}^{\tau}\int_{\mathbb{R}^{4}}\phi_{R,t}(\tau,\xi)\phi_{R,t}(\tau,\tilde{\xi})\sigma'(u(\tau,\xi))\sigma'(u(\tau,\tilde{\xi}))D_{s,y}u(\tau,\xi)D_{s,\tilde{y}}u(\tau,\tilde{\xi})\\
	&\qquad\times|\xi-\tilde{\xi}|^{-\beta}|y-\tilde{y}|^{-\beta} \,d\xi \,d\tilde{\xi}\,dy\,d\tilde{y}\, ds\,d\tau\bigg]\\
	\le&\frac{C}{\sigma_{R,t}^2} \int_{0}^{t}\int_{0}^{\tau}\int_{\mathbb{R}^{4}}\phi_{R,t}(\tau,\xi)\phi_{R,t}(\tau,\tilde{\xi})\mathbb{E}\left[\left|\sigma'(u(\tau,\xi))\sigma'(u(\tau,\tilde{\xi}))D_{s,y}u(\tau,\xi)D_{s,\tilde{y}}u(\tau,\tilde{\xi})\right|\right]\\
	&\qquad\times|\xi-\tilde{\xi}|^{-\beta}|y-\tilde{y}|^{-\beta}\,d\xi \,d\tilde{\xi}\,dy\,d\tilde{y}\, ds\,d\tau\,. 
\end{align*}
Using the H\"{o}lder inequality, the estimate \eqref{ub1} and the assumption \textbf{(H)}, we have
\begin{equation*}
	\mathbb{E}\left[\left|\sigma'(u(\tau,\xi))\sigma'(u(\tau,\tilde{\xi}))D_{s,y}u(\tau,\xi)D_{s,\tilde{y}}u(\tau,\tilde{\xi})\right|\right]\le CG(\tau-s,\xi-y)G(\tau-s,\tilde{\xi}-\tilde{y})\,.
\end{equation*}
Then using the estimate \eqref{Phi} yields
\begin{equation}
	J_2\le \begin{cases}
		c_{t,\beta}R^{-\beta}\,,&\beta\in(0,\frac{1}{2})\,,\\
		c_{t,\beta}R^{\beta-1}\log R\,,&\beta=\frac{1}{2}\,,\\
		c_{t,\beta}R^{\beta-1}\,,&\beta\in(\frac{1}{2},1)\,.
	\end{cases}
\end{equation}
The computations of $J_1$ and $J_2$ yields that
\[
h(R)=\mathbb{E}\left[\|D_{s,y}F_{R,t}\|_{\mathcal{H}}^2\right]
=\begin{cases}
		1+c_{t,\beta}R^{-\beta}\,,&\beta\in(0,\frac{1}{2})\,,\\
		1+c_{t,\beta}R^{\beta-1}\log R\,,&\beta=\frac{1}{2}\,,\\
		1+c_{t,\beta}R^{\beta-1}\,,&\beta\in(\frac{1}{2},1)\,.
	\end{cases}
\] 
So 
\begin{align*}
	\left\|\delta(D_{s,y}F_{R,t})-F_{R,t}\mathbb{E}\left[\|D_{s,y}F_{R,t}\|_{\mathcal{H}}^2\right]\right\|_2^2 \le&h(R)-h(R)^2+\mathbb{E}\left[\| D_{r,z}D_{s,y}F_{R,t}\|_{\mathcal{H}^{\otimes 2}}\right]\\
	\le&\begin{cases}
		c_{t,\beta}R^{-\beta}\,,&\beta\in(0,\frac{1}{2})\,,\\
		c_{t,\beta}R^{\beta-1}\log R\,,&\beta=\frac{1}{2}\,,\\
		c_{t,\beta}R^{\beta-1}\,,&\beta\in(\frac{1}{2},1)\,.
	\end{cases}
\end{align*}
Finally, we get 
\begin{equation}
	\left\|\delta(D_{s,y}F_{R,t})-F_{R,t}\mathbb{E}\left[\|D_{s,y}F_{R,t}\|_{\mathcal{H}}^2\right]\right\|_2\le \begin{cases}
		c_{t,\beta}R^{-\beta/2}\,,&\beta\in(0,\frac{1}{2})\,,\\
		c_{t,\beta}R^{(\beta-1)/2}(\log R)^{1/2}\,,&\beta=\frac{1}{2}\,,\\
		c_{t,\beta}R^{(\beta-1)/2}\,,&\beta\in(\frac{1}{2},1)\,.
	\end{cases}\label{e.3.17} 
\end{equation}
\begin{remark}\label{r.3.3} 
    In the case of space-time white noise, the above term $J_2$ does not decay to zero as $R \to \infty$. This can be understood from two perspectives. On  one hand, as observed in \cite[page 14]{Nualart_Zheng_SWEfractional_2020}, for space white noise, even in the linear case, the first chaos is not the only component contributing to the limiting variance.
On the other hand, a direct calculation yields
\begin{equation*}
\int_{0}^{t}\int_{0}^{\tau}\int_{\mathbb{R}^{2}}\phi_{R,t}^2(\tau,\xi)G^2(\tau-s,\xi-y)
	 \,dy\,d\xi \,ds\,d\tau\sim \frac{4}{15}t^5R\, .
	\end{equation*}
    Moreover, $\sigma_{R,t}^2\sim R$ as $R\to\infty$ (see \cite[Proposition 3.2]{Nualart_Zheng_SWEfractional_2020}). Therefore, after normalization by $\sigma_{R,t}^2$, the corresponding quantity remains of constant order and does not converge to zero.
Consequently, in the case of space-time white noise, the above argument  does not provide the decay estimate required for the term considered there and other approach is needed. 
\end{remark}
\subsection{Estimation of $\sqrt{\mathrm{Var}\left(\|D_{s,y}F_{R,t}\|_{\mathcal{H}}^2 \right)}$}\label{estimation of part 2}
We write
\begin{equation*}
\|D_{s,y}F_{R,t}\|_{\mathcal{H}}^2:=B_1+B_2+B_3\,,
\end{equation*}
where
\begin{equation*}
	B_1=\frac{1}{\sigma_{R,t}^2}\int_{0}^{t}\int_{\mathbb{R}^{2}}\phi_{R,t}(s,y)\phi_{R,t}(s,\tilde{y})\sigma(u(s,y))\sigma(u(s,\tilde{y}))|y-\tilde{y}|^{-\beta}\,dy\,d\tilde{y}\,ds\,,
\end{equation*}
\begin{equation*}
	\begin{split}
		B_2=
		&\frac{2}{\sigma_{R,t}^2}\int_{0}^{t}\int_{\mathbb{R}^{2}}\left(\int_{[s,t]\times\mathbb{R}}\phi_{R,t}(\tau,\xi)\sigma'(u(\tau,\xi))D_{s,y}u(\tau,\xi)\,W(d\tau,d\xi)\right)\\
		&\qquad\times \phi_{R,t}(s,\tilde{y})\sigma(u(s,\tilde{y}))|y-\tilde{y}|^{-\beta}\,dy\,d\tilde{y}\,ds
	\end{split} 
\end{equation*}
and 
\begin{align*}
	B_3=&\frac{1}{\sigma_{R,t}^2}\int_{0}^{t}\int_{\mathbb{R}^{2}}\left(\int_{[s,t]\times\mathbb{R}}\phi_{R,t}(\tau,\xi)\sigma'(u(\tau,\xi))D_{s,y}u(\tau,\xi)\,W(d\tau,d\xi)\right)\\
	&\qquad\times\left(\int_{[s,t]\times\mathbb{R}}\phi_{R,t}(\tau,\xi)\sigma'(u(\tau,\xi))D_{s,\tilde{y}}u(\tau,\xi)\,W(d\tau,d\xi)\right)|y-\tilde{y}|^{-\beta}\,dy\,d\tilde{y}\,ds\,.
\end{align*}
First,  we have 
\begin{equation*}
	\sqrt{\mathrm{Var}\left(\|D_{s,y}F_{R,t}\|_{\mathcal{H}}^2 \right)}\le\sqrt{3}\left(\sqrt{\mathrm{Var}(B_1) }+\sqrt{\mathrm{Var}(B_2) }+\sqrt{\mathrm{Var}(B_3) }\right)\le\sqrt{3}(A_1+A_2+A_3)\,,
\end{equation*}
with
\begin{align*}
	A_1=&\frac{1}{\sigma_{R,t}^2}\int_{0}^{t}\bigg(\int_{\mathbb{R}^{4}}\phi_{R,t}(s,y)\phi_{R,t}(s,\tilde{y})\phi_{R,t}(s,y')\phi_{R,t}(s,\tilde{y}')|y-\tilde{y}|^{-\beta}|y'-\tilde{y}'|^{-\beta}\\
	&\qquad\times\mathrm{Cov}\left[\sigma(u(s,y))\sigma(u(s,\tilde{y})),\sigma(u(s,y'))\sigma(u(s,\tilde{y}'))\right]\,dy\,d\tilde{y}\,dy'\,d\tilde{y}'\bigg)^{1/2}\,ds\,,
\end{align*}
\begin{align*}
	A_2=&\frac{2}{\sigma_{R,t}^2}\int_{0}^{t}\bigg(\int_{\mathbb{R}^{6}}\int_{s}^{t}\phi_{R,t}(\tau,\xi)\phi_{R,t}(\tau,\tilde{\xi})\phi_{R,t}(s,\tilde{y})\phi_{R,t}(s,\tilde{y}')\\
	&\qquad\times\mathbb{E}\left[\sigma'(u(\tau,\xi))D_{s,y}u(\tau,\xi)\sigma'(u(\tau,\tilde{\xi}))D_{s,y'}u(\tau,\tilde{\xi})\sigma(u(s,\tilde{y}))\sigma(u(s,\tilde{y}'))\right]\\
	&\qquad\quad\times|y-\tilde{y}|^{-\beta}|y'-\tilde{y}'|^{-\beta}|\xi-\tilde{\xi}|^{-\beta}\,dy\,d\tilde{y}\,dy'\,d\tilde{y}'\,d\xi \,d\tilde{\xi}\,d\tau\bigg)^{1/2}\,ds
\end{align*}
and
\begin{align*}
	A_3=&\frac{1}{\sigma_{R,t}^2}\int_{0}^{t}\bigg(\int_{\mathbb{R}^{8}}\int_{s}^{t}\int_{s}^{t}\phi_{R,t}(\tau,\xi)\phi_{R,t}(\tau,\tilde{\xi})\phi_{R,t}(\tau',\xi')\phi_{R,t}(\tau',\tilde{\xi}')\\
	&\qquad\times\mathrm{Cov}\left[\sigma'(u(\tau,\xi))D_{s,y}u(\tau,\xi)\sigma'(u(\tau,\tilde{\xi}))D_{s,\tilde{y}}u(\tau,\tilde{\xi}),\sigma'(u(\tau',\xi'))\right. \\
	&\qquad \qquad\left. \times  D_{s,y'}u(\tau',\xi')\sigma'(u(\tau',\tilde{\xi}'))D_{s,\tilde{y}'}u(\tau',\tilde{\xi}')\right]\\
	&\qquad \times|y-\tilde{y}|^{-\beta}|y'-\tilde{y}'|^{-\beta}|\xi-\tilde{\xi}|^{-\beta}|\xi'-\tilde{\xi}'|^{-\beta}\,dy\,d\tilde{y}\,dy'\,d\tilde{y}'\,d\xi \,d\tilde{\xi}\,d\xi' \,d\tilde{\xi}'\,d\tau \,d\tau' \bigg)^{1/2}ds\,.
\end{align*}
For $A_1$ and $A_2$, we can obtain the following estimates similar to that in the proof of   \cite[Theorem 1.1 ]{Nualart_Zheng_SWEfractional_2020} \begin{equation}
	A_1\le CR^{-\frac{\beta}{2}},\quad  A_2\le CR^{-\frac{\beta}{2}}\quad    \mathrm{for}\ R\ \mathrm{big} \ \mathrm{enough}\,.
\end{equation}
As for $A_3$, we can write
\begin{align*}
	&\left|\mathrm{Cov}\left[\sigma'(u(\tau,\xi))D_{s,y}u(\tau,\xi)\sigma'(u(\tau,\tilde{\xi}))D_{s,\tilde{y}}u(\tau,\tilde{\xi}),\sigma'(u(\tau',\xi'))D_{s,y'}u(\tau',\xi')\sigma'(u(\tau',\tilde{\xi}'))D_{s,\tilde{y}'}u(\tau',\tilde{\xi}')\right]\right|\\
	=&\bigg|\mathbb{E}\left[\sigma'(u(\tau,\xi))D_{s,y}u(\tau,\xi)\sigma'(u(\tau,\tilde{\xi}))D_{s,\tilde{y}}u(\tau,\tilde{\xi})\sigma'(u(\tau',\xi'))D_{s,y'}u(\tau',\xi')\sigma'(u(\tau',\tilde{\xi}'))D_{s,\tilde{y}'}u(\tau',\tilde{\xi}') \right]\\
	&-\mathbb{E}\left[\sigma'(u(\tau,\xi))D_{s,y}u(\tau,\xi)\sigma'(u(\tau,\tilde{\xi}))D_{s,\tilde{y}}u(\tau,\tilde{\xi})\right]\mathbb{E}\left[\sigma'(u(\tau',\xi'))D_{s,y'}u(\tau',\xi')\sigma'(u(\tau',\tilde{\xi}'))D_{s,\tilde{y}'}u(\tau',\tilde{\xi}')\right]\bigg|\\
	\le&\left|\mathbb{E}\left[\sigma'(u(\tau,\xi))D_{s,y}u(\tau,\xi)\sigma'(u(\tau,\tilde{\xi}))D_{s,\tilde{y}}u(\tau,\tilde{\xi})\sigma'(u(\tau',\xi'))D_{s,y'}u(\tau',\xi')\sigma'(u(\tau',\tilde{\xi}'))D_{s,\tilde{y}'}u(\tau',\tilde{\xi}') \right]\right|\\
	&+\left|\mathbb{E}\left[\sigma'(u(\tau,\xi))D_{s,y}u(\tau,\xi)\sigma'(u(\tau,\tilde{\xi}))D_{s,\tilde{y}}u(\tau,\tilde{\xi})\right]\right|\left|\mathbb{E}\left[\sigma'(u(\tau',\xi'))D_{s,y'}u(\tau',\xi')\sigma'(u(\tau',\tilde{\xi}'))D_{s,\tilde{y}'}u(\tau',\tilde{\xi}')\right]\right|\,.
\end{align*}
Using H\"{o}lder's inequality, assumption \textbf{(H)} and the estimate \eqref{ub1}, we have
\begin{align*}
	&\left|\mathbb{E}\left[\sigma'(u(\tau,\xi))D_{s,y}u(\tau,\xi)\sigma'(u(\tau,\tilde{\xi}))D_{s,\tilde{y}}u(\tau,\tilde{\xi})\sigma'(u(\tau',\xi'))D_{s,y'}u(\tau',\xi')\sigma'(u(\tau',\tilde{\xi}'))D_{s,\tilde{y}'}u(\tau',\tilde{\xi}') \right]\right|\\
	\le&CG(\tau-s,\xi-y)G(\tau-s,\tilde{\xi}-\tilde{y})G(\tau'-s,\xi'-y')G(\tau'-s,\tilde{\xi}'-\tilde{y}') \,.
\end{align*}
Similarly,
\begin{align*}
	\left|\mathbb{E}\left[\sigma'(u(\tau,\xi))D_{s,y}u(\tau,\xi)\sigma'(u(\tau,\tilde{\xi}))D_{s,\tilde{y}}u(\tau,\tilde{\xi})\right]\right|\le C G(\tau-s,\xi-y)G(\tau-s,\tilde{\xi}-\tilde{y})
\end{align*}
and
\begin{align*}
	&\left|\mathbb{E}\left[\sigma'(u(\tau',\xi'))D_{s,y'}u(\tau',\xi')\sigma'(u(\tau',\tilde{\xi}'))D_{s,\tilde{y}'}u(\tau',\tilde{\xi}')\right]\right|\\
	\le& C G(\tau'-s,\xi'-y')G(\tau'-s,\tilde{\xi}'-\tilde{y}')\,.
\end{align*}
Then applying the same methods used for the estimate \eqref{Phi}, we obtain
\begin{align*}
	A_3=&\frac{C}{\sigma_{R,t}^2}\int_{0}^{t}\bigg(\int_{s}^{t}\int_{s}^{t}\int_{\mathbb{R}^{8}}\phi_{R,t}(\tau,\xi)\phi_{R,t}(\tau,\tilde{\xi})\phi_{R,t}(\tau',\xi')\phi_{R,t}(\tau',\tilde{\xi}')G(\tau-s,\xi-y)\\
	&\qquad\times G(\tau-s,\tilde{\xi}-\tilde{y})G(\tau'-s,\xi'-y')G(\tau'-s,\tilde{\xi}'-\tilde{y}') |y-\tilde{y}|^{-\beta}|y'-\tilde{y}'|^{-\beta}|\xi-\tilde{\xi}|^{-\beta}\\
	&\qquad\quad\times|\xi'-\tilde{\xi}'|^{-\beta}\,dy\,d\tilde{y}\,dy'\,d\tilde{y}'\,d\xi \,d\tilde{\xi}\,d\xi' \,d\tilde{\xi}'\,d\tau \,d\tau' \bigg)^{1/2}ds\\
	\le&\begin{cases}
		c_{t,\beta}R^{-\beta}\,,&\beta\in(0,\frac{1}{2})\,,\\
		c_{t,\beta}R^{\beta-1}\log R\,,&\beta=\frac{1}{2}\,,\\
		c_{t,\beta}R^{\beta-1}\,,&\beta\in(\frac{1}{2},1)\,.
	\end{cases}
\end{align*}\
Summarizing the above argument we have
\begin{equation}
	\sqrt{\mathrm{Var}\left(\|D_{s,y}F_{R,t}\|_{\mathcal{H}}^2 \right)}
	\le \begin{cases}
		c_{t,\beta}R^{-\beta/2}\,,&\beta\in(0,\frac{2}{3}]\,,\\
		c_{t,\beta}R^{\beta-1}\,,&\beta\in(\frac{2}{3},1)\,.
	\end{cases} \label{e.3.19}
\end{equation}

\subsection{Estimation of $\sqrt{ \mathrm{Var}[\langle D_{s,y}F_{R,t},w\rangle_{\mathcal{H}}]}$}\label{estimation of part 4}
The following quantitative central limit theorem has been proved in \cite[Theorem 1.1]{Nualart_Zheng_SWEfractional_2020}.
\begin{theorem}
	Let $d_{TV}$ denote the total variation distance and let $Z\sim N(0,1)$. Let $F_{R,t}=\delta(w)$ and for any fixed $t>0$, there exists a constant $C=C_{t,\beta,\sigma}$, depending on $t,\beta$ and $\sigma$, such that
	\begin{equation}
		d_{TV}(F_{R,t},Z)\le 2\sqrt{ \mathrm{Var}[\langle D_{s,y}F_{R,t},w\rangle_{\mathcal{H}}]}\le CR^{-\frac{\beta}{2}}\,.
		\label{e.3.21}
	\end{equation}
\end{theorem}

\subsection{Proof of Theorem \ref{main result}}\label{final proof}
We are going to use \eqref{ub of density}. The first factor in 
\eqref{ub of density} is bounded by \eqref{e.3.14}  with $q=1$.  
The first, second, and third
term inside the bracket are bounded by
\eqref{e.3.17}, \eqref{e.3.19},
and \eqref{second detivative estimate}. The last term in 
\eqref{ub of density} is bounded by \eqref{e.3.21}.  
This proves Theorem \ref{main result}.

\section{Appendix}\label{appendix}
\begin{lemma}
		Let $\Phi_{R,t}$ be as in \eqref{Phi1}. For $0<r<s<t$ and $x,y,z\in\mathbb{R}$,
	we have 
\begin{align}\label{Phi}
	\Phi_{R,t}\le\begin{cases}
		c_{t,\beta}R^{-\beta}\,,&\beta\in(0,\frac{1}{2})\,,\\
		c_{t,\beta}R^{\beta-1}\log R\,,&\beta=\frac{1}{2}\,,\\
		c_{t,\beta}R^{\beta-1}\,,&\beta\in(\frac{1}{2},1)\,.
	\end{cases}
\end{align}
\end{lemma}
\begin{proof}
We can write
\begin{align*}
	\Phi_{R,t}
	= &\frac{1}{\sigma_{R,t}^2}\int_0^t \int_0^s \int_{\mathbb{R}^4} \phi_{R,t}(s,y)\phi_{R,t}(s,y')
	G(s-r,y-z)G(s-r,y'-z') \\
	&\qquad \times |y-y'|^{-\beta}\,|z-z'|^{-\beta}\,dz\,dz'\,dy\,dy'\,dr\,ds \\
	= &\frac{C}{\sigma_{R,t}^2} \int_0^t \int_0^s \int_{\mathbb{R}^4}\int_{[-R,R]^2}
	\mathbf{1}_{\{|x-y|\le t-s\}}\mathbf{1}_{\{|x'-y'|\le t-s\}}
	\mathbf{1}_{\{|y-z|\le s-r\}}\mathbf{1}_{\{|y'-z'|\le s-r\}} \\
	&\qquad \times |y-y'|^{-\beta}\,|z-z'|^{-\beta}\,dx\,dx'\,dz\,dz'\,dy\,dy'\,dr\,ds\,.
\end{align*}
Assuming that $R \ge t$ and integrating in the variables $x,x' \in [-R,R]$, we have
\begin{align*}
	\Phi_{R,t}
	\le& \frac{C_t}{\sigma_{R,t}^2} \int_0^t \int_0^s \int_{\mathbb{R}^4}
	\mathbf{1}_{\{|y-z|\le s-r,\,|y'-z'|\le s-r\}}
	\mathbf{1}_{[-2R,2R]}(y)\mathbf{1}_{[-2R,2R]}(y') \\
	&\qquad \times |y-y'|^{-\beta}\,|z-z'|^{-\beta}\,dz\,dz'\,dy\,dy'\,dr\,ds \\
	= &\frac{C_{t,\beta}}{\sigma_{R,t}^2} \int_0^t \int_0^s \int_{[-2R,2R]^2}
	\Big( |y-y'-2(s-r)|^{2-\beta} + |y-y'+2(s-r)|^{2-\beta}
	- 2|y-y'|^{2-\beta}\Big) \\
	&\qquad \times |y-y'|^{-\beta}\,dy\,dy'\,dr\,ds\,,
\end{align*}
where the last equality follows from Lemma \ref{int1}.
Denote that $\Delta=y-y', b=s-r>0$, and
\begin{align*}
	B_{\beta}(\Delta;b):=&\left|\Delta-2b\right|^{2-\beta}+\left|\Delta+2b\right|^{2-\beta}-2|\Delta|^{2-\beta}\\
	=&(2-\beta)\int_{0}^{2b}\left(\mathrm{sgn}(\Delta+\theta)|\Delta+\theta|^{1-\beta}-\mathrm{sgn}(\Delta-\theta)|\Delta-\theta|^{1-\beta}\right)\,d\theta\\
	=&(2-\beta)(1-\beta)\int_{0}^{2b}\int_{\Delta-\theta}^{\Delta+\theta}|\eta|^{-\beta}\,d\eta \,d\theta\\
	=&(2-\beta)(1-\beta)\int_{0}^{2b}\int_{-\theta}^{\theta}|\Delta+\eta|^{-\beta}\,d\eta \,d\theta\\
	=&(2-\beta)(1-\beta)\int_{-2b}^{2b}\left(\int_{|\eta|}^{2b}d\theta\right)|\Delta+\eta|^{-\beta}\,d\eta\\
	=&(2-\beta)(1-\beta)\int_{-2b}^{2b}(2b-|\eta|)|\Delta+\eta|^{-\beta}\,d\eta\,.
\end{align*}
Clearly, $	B_{\beta}(\Delta;b)\ge0 $. Then we will show that
		\begin{align}
			B_{\beta}(\Delta;b)\le C_{\beta}\begin{cases}
			    b^{2-\beta} & \hbox{when $|\Delta|< 4b$}\\
                b^2|\Delta|^{-\beta} & \hbox{when $|\Delta|\ge  4b$}\,. \\  
			\end{cases} \label{e.4.55} 
		\end{align}
		When $|\Delta|\ge 4b$, since $|\eta|\le2b$, we have $|\Delta\pm \eta|\ge\frac{1}{2}|\Delta|$, then
		\begin{align*}
			B_{\beta}(\Delta;b)=&(2-\beta)(1-\beta)\int_{-2b}^{2b}(2b-|\eta|)|\Delta+\eta|^{-\beta}\,d\eta\\
			\le &C_{\beta}b^2|\Delta|^{-\beta}\,.
		\end{align*}
		When $ |\Delta|<4b$, we have
		\begin{align*}
			B_{\beta}(\Delta;b)=&\left|\Delta-2b\right|^{2-\beta}+\left|\Delta+2b\right|^{2-\beta}-2|\Delta|^{2-\beta}\\ \le&2(|\Delta|+2b)^{2-\beta}-2|\Delta|^{2-\beta}\\
			\le&2(|\Delta|+2b)^{2-\beta}\\
			\le& C_{\beta}b^{2-\beta}\,.
		\end{align*}
		Now we have
		\begin{align*}
			\Phi_{R,t}\le \frac{C_{t,\beta}}{\sigma_{R,t}^2} \int_0^t \int_0^s \int_{[-2R,2R]^2}
			|y-y'|^{-\beta}	B_{\beta}(y-y';b)\,dy\,dy'\,dr\,ds\,.
		\end{align*}
		Let $u=y-y',v=(y+y')/2$. The integration region becomes $\{(u,v):|u|\le 4R,|v|\le2R-|u|/2\}$. Therefore,
		\begin{equation*}
			\int_{[-2R,2R]^2}
			|y-y'|^{-\beta}B_{\beta}(y-y';b)\,dy\,dy'=\int_{-4R}^{4R}(4R-|u|)|u|^{-\beta}B_{\beta}(u;b)\,du\,.
		\end{equation*}
		Then we obtain
		\begin{align*}
			\Phi_{R,t}\le &\frac{C_{t,\beta}}{\sigma_{R,t}^2}\int_0^t \int_0^s \int_{-4R}^{4R}(4R-|u|)|u|^{-\beta}B_{\beta}(\Delta;b) \,du\,dr\,ds\\
			\le &\frac{C_{t,\beta}}{\sigma_{R,t}^2}\int_{0}^{t}\int_{0}^{s}\left[b^{2-\beta}\int_{|u|\le 4b}(4R-|u|)|u|^{-\beta}\,du+b^2\int_{4b<|u|\le 4R}(4R-|u|)|u|^{-2\beta}\,du\right]\,dr\,ds\\
			=&:\frac{C_{t,\beta}}{\sigma_{R,t}^2}\int_{0}^{t}\int_{0}^{s}(I_1(b,R)+I_2(b,R))\,dr\,ds\,.
		\end{align*}
		For $I_1(b,R) $, directly,
		\begin{equation}\label{re1}
			I_1(b,R)\le 4Rb^{2-\beta}\int_{0}^{4b}u^{-\beta}\,du\le c_{t,\beta}R\,.
		\end{equation}
		As for $I_2(b,R) $,
		\begin{align}\label{re2}
			I_2(b,R)\le 4Rb^2\int_{4b}^{4R}u^{-2\beta}\,du\le
			\begin{cases}
				c_{t,\beta}R^{2-2\beta},&\beta\in(0,\frac{1}{2}),\\
				c_{t,\beta}R\log R,&\beta=\frac{1}{2},\\
				c_{t,\beta}R,&\beta\in(\frac{1}{2},1)\,.
			\end{cases}
		\end{align}
		Finally, using (20), we have
		\begin{align*}
			\Phi_{R,t}\le\begin{cases}
				c_{t,\beta}R^{-\beta}\,,&\beta\in(0,\frac{1}{2})\,,\\
				c_{t,\beta}R^{\beta-1}\log R\,,&\beta=\frac{1}{2}\,,\\
				c_{t,\beta}R^{\beta-1}\,,&\beta\in(\frac{1}{2},1)\,.
			\end{cases}
		\end{align*}
        This proves the lemma. 
\end{proof}
\begin{lemma}
		Let $0<r<s<\tau<t$ and $x,y,\xi,z\in\mathbb{R}$.\\
	(i) Let $\Phi_{R,t}^{(1)}$ be as in \eqref{Phi2}.
We have 
\begin{align}\label{Phi_2}
	\Phi_{R,t}^{(1)}\le\begin{cases}
		c_{t,\beta}R^{-\beta}\,,&\beta\in(0,\frac{1}{2})\,,\\
		c_{t,\beta}R^{\beta-1}\log R\,,&\beta=\frac{1}{2}\,,\\
		c_{t,\beta}R^{\beta-1}\,,&\beta\in(\frac{1}{2},1)\,.
	\end{cases}
\end{align}
(ii) Let $\Phi_{R,t}^{(2)}$ be as in \eqref{Phi3}.
We have 
\begin{align}\label{Phi_3}
	\Phi_{R,t}^{(2)}\le\begin{cases}
		c_{t,\beta}R^{-\beta}\,,&\beta\in(0,\frac{1}{2})\,,\\
		c_{t,\beta}R^{\beta-1}\log R\,,&\beta=\frac{1}{2}\,,\\
		c_{t,\beta}R^{\beta-1}\,,&\beta\in(\frac{1}{2},1)\,.
	\end{cases}
\end{align}
\end{lemma}
\begin{proof}
We can write
\begin{align*}
	\Phi_{R,t}^{(1)}=&\frac{1}{\sigma_{R,t}^2}\int_{0}^{t}\int_{0}^{s}\int_{s}^{t}\int_{\mathbb{R}^6}\phi_{R,t}(\tau,\xi)\phi_{R,t}(\tau,\xi')G(\tau-r,\xi-z)G(\tau-r,\xi'-z')G(\tau-s,\xi-y)\\
	&\times G(\tau-s,\xi'-y')|\xi-\xi'|^{-\beta}|z-z'|^{-\beta}|y-y'|^{-\beta}\,d\xi\, d\xi'\,dz\,dz'\,dy\,dy'\,d\tau \,dr\,ds\\
	=&\frac{1}{\sigma_{R,t}^2}\int_{0}^{t}\int_{0}^{s}\int_{s}^{t}\int_{\mathbb{R}^6}\int_{[-R,R]^2}G(t-\tau,x_1-\xi)G(t-\tau,x_2-\xi')G(\tau-r,\xi-z)G(\tau-r,\xi'-z')\\
&\times G(\tau-s,\xi-y)G(\tau-s,\xi'-y')|\xi-\xi'|^{-\beta}|z-z'|^{-\beta}|y-y'|^{-\beta}\,dx_1\,dx_2\,d\xi\, d\xi'\,dz\,dz'\,dy\,dy'\,d\tau \,dr\,ds\\
=&\frac{C}{\sigma_{R,t}^2}\int_{0}^{t}\int_{0}^{s}\int_{s}^{t}\int_{\mathbb{R}^6}\int_{[-R,R]^2}\mathbf{1}_{\{|x_1-\xi|<t-\tau\}}\mathbf{1}_{\{|x_2-\xi'|<t-\tau\}}\mathbf{1}_{\{|\xi-y|<\tau-s\}}\mathbf{1}_{\{|\xi'-y'|<\tau-s\}}\\
&\times \mathbf{1}_{\{|\xi-z|<\tau-r\}}\mathbf{1}_{\{|\xi'-z'|<\tau-r\}}|\xi-\xi'|^{-\beta}|z-z'|^{-\beta}|y-y'|^{-\beta}\,dx_1\,dx_2\,d\xi\, d\xi'\,dz\,dz'\,dy\,dy'\,d\tau \,dr\,ds\,.
\end{align*}
 Assuming that $R \ge t$ and integrating in the variables $x,x' \in [-R,R]$, we have
\begin{align*}
	\Phi_{R,t}^{(1)}\le&\frac{C_t}{\sigma_{R,t}^2}\int_{0}^{t}\int_{0}^{s}\int_{s}^{t}\int_{\mathbb{R}^6}\mathbf{1}_{[-2R,2R]}(\xi)\mathbf{1}_{[-2R,2R]}(\xi')\mathbf{1}_{\{|\xi-y|<\tau-s\}}\mathbf{1}_{\{|\xi'-y'|<\tau-s\}}\\
	&\times \mathbf{1}_{\{|\xi-z|<\tau-r\}}\mathbf{1}_{\{|\xi'-z'|<\tau-r\}}|\xi-\xi'|^{-\beta}|z-z'|^{-\beta}|y-y'|^{-\beta}\,d\xi\, d\xi'\,dz\,dz'\,dy\,dy'\,d\tau \,dr\,ds\\
	\le&\frac{C_{t,\beta}}{\sigma_{R,t}^2}\int_{0}^{t}\int_{0}^{s}\int_{s}^{t}\int_{\mathbb{R}^4}\mathbf{1}_{[-2R,2R]}(\xi)\mathbf{1}_{[-2R,2R]}(\xi')\mathbf{1}_{\{|\xi-y|<\tau-r\}}\mathbf{1}_{\{|\xi'-y'|<\tau-r\}}\\
	&\times 	\Big( |\xi-\xi'-2(\tau-r)|^{2-\beta} + |\xi-\xi'+2(\tau-r)|^{2-\beta}
	- 2|\xi-\xi'|^{2-\beta}\Big)\\
	&\quad\times|\xi-\xi'|^{-\beta}|y-y'|^{-\beta}\,d\xi\, d\xi'\,dy\,dy'\,d\tau \,dr\,ds\\
	=&\frac{C_{t,\beta}}{\sigma_{R,t}^2}\int_{0}^{t}\int_{0}^{s}\int_{s}^{t}\int_{[-2R,2R]^2}\Big( |\xi-\xi'-2(\tau-r)|^{2-\beta} + |\xi-\xi'+2(\tau-r)|^{2-\beta}
	- 2|\xi-\xi'|^{2-\beta}\Big)^2\\
	&\times |\xi-\xi'|^{-\beta}\,d\xi \,d\xi'\,d\tau \,dr\,ds\,,
\end{align*}
where the last two steps follows from Lemma \ref{int1}.
As for  \eqref{e.4.55}, we have easily 
\begin{align*}
	|\xi-\xi'-2(\tau-r)|^{2-\beta} &+ |\xi-\xi'+2(\tau-r)|^{2-\beta}
	- 2|\xi-\xi'|^{2-\beta}\\
    &\le C_{\beta}
    \begin{cases} (\tau-r)^{2-\beta} &\hbox{when $\xi-\xi'<4(\tau-r)$} \\
    (\tau-r)^2|\xi-\xi'|^{-\beta}&\hbox{when $\xi-\xi'\ge 4(\tau-r)$}  \,.
    \end{cases} 
\end{align*}
Let $\tilde{u}=\xi-\xi',\tilde{v}=(\xi+\xi')/2$. The integration region becomes $\{(\tilde{u},\tilde{v}):|\tilde{u}|\le 4R,|\tilde{v}|\le2R-|\tilde{u}|/2\}$. Therefore,
\begin{equation*}
	\int_{[-2R,2R]^2}
	|\xi-\xi'|^{-\beta}\,d\xi\, d\xi'=\int_{-4R}^{4R}(4R-|\tilde{u}|)|\tilde{u}|^{-\beta}\,d\tilde{u}\,.
\end{equation*}
Let $\tilde{b}=\tau-r$, we obtain
\begin{align*}
	\Phi_{R,t}\le &\frac{C_{t,\beta}}{\sigma_{R,t}^2} \int_0^t \int_0^s\int_{s}^{t} \int_{-4R}^{4R}(4R-|\tilde{u}|)|\tilde{u}|^{-\beta}\min\{\tilde{b}^{4-2\beta},\tilde{b}^{4-\beta}|\tilde{u}|^{-\beta}\}\,du\,d\tau \,dr\,ds\\
	\le &\frac{C_{t,\beta}}{\sigma_{R,t}^2}\int_{0}^{t}\int_{0}^{s}\int_{s}^{t}\left[\tilde{b}^{4-2\beta}\int_{|u|\le 4\tilde{b}}(4R-|\tilde{u}|)|\tilde{u}|^{-\beta}\,d\tilde{u}+\tilde{b}^{4-\beta}\int_{4\tilde{b}<|\tilde{u}|\le4R}(4R-|\tilde{u}|)|\tilde{u}|^{-2\beta}\,d\tilde{u}\right]\,d\tau \,ds\,dr\\
	=&:\frac{C_{t,\beta}}{\sigma_{R,t}^2}\int_{0}^{t}\int_{0}^{s}\int_{s}^{t}(\tilde{I}_1(\tilde{b},R)+\tilde{I}_2(\tilde{b},R))\,d\tau \,dr\,ds\,.
\end{align*}
By a similar argument  we can obtain the   bounds  same as \eqref{re1} and \eqref{re2}
for $\tilde{I}_1(\tilde{b},R)$ and  $\tilde{I}_2(\tilde{b},R)$,  and substituting those estimates  into the above integration
will yield  
\begin{align*}
	\Phi_{R,t}^{(1)}\le\begin{cases}
		c_tR^{-\beta}\,,&\beta\in(0,\frac{1}{2})\,,\\
		c_tR^{\beta-1}\log R\,,&\beta=\frac{1}{2}\,,\\
		c_tR^{\beta-1}\,,&\beta\in(\frac{1}{2},1)\,.
	\end{cases}
\end{align*}
Similarly, the   estimate for $\Phi_{R,t}^{(2)}$ can be obtained by the same argument.
\end{proof}

\bibliographystyle{plain}  
\bibliography{ref}  

@article{Balan_SPA_2025,
title = {Gaussian fluctuations for the wave equation under rough random perturbations},
journal = {Stochastic Processes and their Applications},
volume = {182},
pages = {104569},
year = {2025},
issn = {0304-4149},
doi = {https://doi.org/10.1016/j.spa.2025.104569},
url = {https://www.sciencedirect.com/science/article/pii/S0304414925000080},
author = {Raluca M. Balan and Jingyu Huang and Xiong Wang and Panqiu Xia and Wangjun Yuan}
}

@article {Dalang,
	AUTHOR = {Dalang, Robert C.},
	TITLE = {Extending the martingale measure stochastic integral with
		applications to spatially homogeneous s.p.d.e.'s},
	JOURNAL = {Electron. J. Probab.},
	FJOURNAL = {Electronic Journal of Probability},
	VOLUME = {4},
	YEAR = {1999},
	PAGES = {no. 6, 29},
	ISSN = {1083-6489},
	MRCLASS = {60H05 (35R60 60G15 60G48 60H15)},
	MRNUMBER = {1684157},
	MRREVIEWER = {Marta\ Sanz Sol\'e},
	DOI = {10.1214/EJP.v4-43},
	URL = {https://doi.org/10.1214/EJP.v4-43},
}

@article {Nualart_Zheng_SWEfractional_2020,
	AUTHOR = {Delgado-Vences, Francisco and Nualart, David and Zheng,
	Guangqu},
	TITLE = {A central limit theorem for the stochastic wave equation with
	fractional noise},
	JOURNAL = {Ann. Inst. Henri Poincar\'e{} Probab. Stat.},
	FJOURNAL = {Annales de l'Institut Henri Poincar\'e{} Probabilit\'es et
	Statistiques},
	VOLUME = {56},
	YEAR = {2020},
	NUMBER = {4},
	PAGES = {3020--3042},
	ISSN = {0246-0203,1778-7017},
	MRCLASS = {60H15 (60F05 60G15 60H07)},
	MRNUMBER = {4164864},
	MRREVIEWER = {Roger\ Pettersson},
	DOI = {10.1214/20-AIHP1069},
	URL = {https://doi.org/10.1214/20-AIHP1069},
}

@book {Nualart_Malliavin_2006,
	AUTHOR = {Nualart, David},
	TITLE = {The {M}alliavin Calculus and Related Topics},
	SERIES = {Probability and its Applications (New York)},
	EDITION = {Second},
	PUBLISHER = {Springer-Verlag, Berlin},
	YEAR = {2006},
	PAGES = {xiv+382},
	ISBN = {978-3-540-28328-7; 3-540-28328-5},
	MRCLASS = {60-02 (60H07 60H30)},
	MRNUMBER = {2200233},
	MRREVIEWER = {Daniel\ Ocone},
}

@article {Balan_Nualart_Quer-Sardanyons_2022,
	AUTHOR = {Balan, Raluca M. and Nualart, David and Quer-Sardanyons,
	Llu\'is and Zheng, Guangqu},
	TITLE = {The hyperbolic {A}nderson model: moment estimates of the
	{M}alliavin derivatives and applications},
	JOURNAL = {Stoch. Partial Differ. Equ. Anal. Comput.},
	FJOURNAL = {Stochastic Partial Differential Equations. Analysis and
	Computations},
	VOLUME = {10},
	YEAR = {2022},
	NUMBER = {3},
	PAGES = {757--827},
	ISSN = {2194-0401,2194-041X},
	MRCLASS = {60H15 (35K57 60F05 60G15 60H07)},
	MRNUMBER = {4491503},
	DOI = {10.1007/s40072-021-00227-5},
	URL = {https://doi.org/10.1007/s40072-021-00227-5},
}

@article {Nualart_density_2007,
	AUTHOR = {Nualart, David and Quer-Sardanyons, Llu\'is},
	TITLE = {Existence and smoothness of the density for spatially
	homogeneous {SPDE}s},
	JOURNAL = {Potential Anal.},
	FJOURNAL = {Potential Analysis. An International Journal Devoted to the
	Interactions between Potential Theory, Probability Theory,
	Geometry and Functional Analysis},
	VOLUME = {27},
	YEAR = {2007},
	NUMBER = {3},
	PAGES = {281--299},
	ISSN = {0926-2601,1572-929X},
	MRCLASS = {60H07 (60H15)},
	MRNUMBER = {2336301},
	MRREVIEWER = {David\ M\'arquez-Carreras},
	DOI = {10.1007/s11118-007-9055-3},
	URL = {https://doi.org/10.1007/s11118-007-9055-3},
}

@article {Huang_Nualart_CLTforSHE_2020,
	AUTHOR = {Huang, Jingyu and Nualart, David and Viitasaari, Lauri},
	TITLE = {A central limit theorem for the stochastic heat equation},
	JOURNAL = {Stochastic Process. Appl.},
	FJOURNAL = {Stochastic Processes and their Applications},
	VOLUME = {130},
	YEAR = {2020},
	NUMBER = {12},
	PAGES = {7170--7184},
	ISSN = {0304-4149,1879-209X},
	MRCLASS = {60H15 (60F05 60G15 60H07)},
	MRNUMBER = {4167203},
	MRREVIEWER = {Zhi\ Wang},
	DOI = {10.1016/j.spa.2020.07.010},
	URL = {https://doi.org/10.1016/j.spa.2020.07.010},
}

@article {Huang_Nualart_SHEcolored_2020,
	AUTHOR = {Huang, Jingyu and Nualart, David and Viitasaari, Lauri and
	Zheng, Guangqu},
	TITLE = {Gaussian fluctuations for the stochastic heat equation with
	colored noise},
	JOURNAL = {Stoch. Partial Differ. Equ. Anal. Comput.},
	FJOURNAL = {Stochastic Partial Differential Equations. Analysis and
	Computations},
	VOLUME = {8},
	YEAR = {2020},
	NUMBER = {2},
	PAGES = {402--421},
	ISSN = {2194-0401,2194-041X},
	MRCLASS = {60H15 (35K91 35R60 60F05 60F17 60G15 60H07)},
	MRNUMBER = {4098872},
	MRREVIEWER = {Stefano\ Bonaccorsi},
	DOI = {10.1007/s40072-019-00149-3},
	URL = {https://doi.org/10.1007/s40072-019-00149-3},
}

@article {Nualart_Xia_Zheng_PAMcolored_2022,
	AUTHOR = {Nualart, David and Xia, Panqiu and Zheng, Guangqu},
	TITLE = {Quantitative central limit theorems for the parabolic
	{A}nderson model driven by colored noises},
	JOURNAL = {Electron. J. Probab.},
	FJOURNAL = {Electronic Journal of Probability},
	VOLUME = {27},
	YEAR = {2022},
	PAGES = {Paper No. 120, 43},
	ISSN = {1083-6489},
	MRCLASS = {60H15 (60F05 60G22 60H07)},
	MRNUMBER = {4479916},
	MRREVIEWER = {Adam\ Matthew\ Bowditch},
	DOI = {10.1214/22-ejp847},
	URL = {https://doi.org/10.1214/22-ejp847},
}

@article {Nualart_Song_Zheng_PAMrough_2020,
	AUTHOR = {Nualart, David and Song, Xiaoming and Zheng, Guangqu},
	TITLE = {Spatial averages for the parabolic {A}nderson model driven by
	rough noise},
	JOURNAL = {ALEA Lat. Am. J. Probab. Math. Stat.},
	FJOURNAL = {ALEA. Latin American Journal of Probability and Mathematical
	Statistics},
	VOLUME = {18},
	YEAR = {2021},
	NUMBER = {1},
	PAGES = {907--943},
	ISSN = {1980-0436},
	MRCLASS = {60H15 (60F05 60G15 60H07)},
	MRNUMBER = {4243520},
	MRREVIEWER = {Yong\ Liu},
	DOI = {10.30757/alea.v18-33},
	URL = {https://doi.org/10.30757/alea.v18-33},
}

@article {Nualart_Zheng_SHEtimefractional_2020,
	AUTHOR = {Nualart, David and Zheng, Guangqu},
	TITLE = {Averaging {G}aussian functionals},
	JOURNAL = {Electron. J. Probab.},
	FJOURNAL = {Electronic Journal of Probability},
	VOLUME = {25},
	YEAR = {2020},
	PAGES = {Paper No. 48, 54},
	ISSN = {1083-6489},
	MRCLASS = {60H15 (60F05 60G15 60H07)},
	MRNUMBER = {4092767},
	MRREVIEWER = {Ciprian\ A.\ Tudor},
	DOI = {10.1214/20-ejp453},
	URL = {https://doi.org/10.1214/20-ejp453},
}

@article {Balan_Yuan_SHEtimeindependent_2023,
	AUTHOR = {Balan, Raluca M. and Yuan, Wangjun},
	TITLE = {Central limit theorems for heat equation with time-independent
	noise: the regular and rough cases},
	JOURNAL = {Infin. Dimens. Anal. Quantum Probab. Relat. Top.},
	FJOURNAL = {Infinite Dimensional Analysis, Quantum Probability and Related
	Topics},
	VOLUME = {26},
	YEAR = {2023},
	NUMBER = {2},
	PAGES = {Paper No. 2250029, 47},
	ISSN = {0219-0257,1793-6306},
	MRCLASS = {60H15 (60F05 60G15 60H07)},
	MRNUMBER = {4598108},
	MRREVIEWER = {Ciprian\ A.\ Tudor},
	DOI = {10.1142/S0219025722500291},
	URL = {https://doi.org/10.1142/S0219025722500291},
}

@article {Nualart_Zheng_2dSWE_2021,
	AUTHOR = {Bola\~nos Guerrero, Raul and Nualart, David and Zheng,
	Guangqu},
	TITLE = {Averaging 2d stochastic wave equation},
	JOURNAL = {Electron. J. Probab.},
	FJOURNAL = {Electronic Journal of Probability},
	VOLUME = {26},
	YEAR = {2021},
	PAGES = {Paper No. 102, 32},
	ISSN = {1083-6489},
	MRCLASS = {60H15 (60F05 60G15 60H07)},
	MRNUMBER = {4290504},
	DOI = {10.1214/21-ejp672},
	URL = {https://doi.org/10.1214/21-ejp672},
}

@article {Nualart_Zheng_SWEd12_2022,
	AUTHOR = {Nualart, David and Zheng, Guangqu},
	TITLE = {Central limit theorems for stochastic wave equations in
	dimensions one and two},
	JOURNAL = {Stoch. Partial Differ. Equ. Anal. Comput.},
	FJOURNAL = {Stochastic Partial Differential Equations. Analysis and
	Computations},
	VOLUME = {10},
	YEAR = {2022},
	NUMBER = {2},
	PAGES = {392--418},
	ISSN = {2194-0401,2194-041X},
	MRCLASS = {60H15 (60F05 60G15 60H07)},
	MRNUMBER = {4439987},
	MRREVIEWER = {Sergey\ V.\ Lototsky},
	DOI = {10.1007/s40072-021-00209-7},
	URL = {https://doi.org/10.1007/s40072-021-00209-7},
}

@article {Balan_Yuan_HAMtimeindependent_2022,
	AUTHOR = {Balan, Raluca M. and Yuan, Wangjun},
	TITLE = {Spatial integral of the solution to hyperbolic {A}nderson
	model with time-independent noise},
	JOURNAL = {Stochastic Process. Appl.},
	FJOURNAL = {Stochastic Processes and their Applications},
	VOLUME = {152},
	YEAR = {2022},
	PAGES = {177--207},
	ISSN = {0304-4149,1879-209X},
	MRCLASS = {60H15 (60F05 60G15 60H07)},
	MRNUMBER = {4450464},
	MRREVIEWER = {Ciprian\ A.\ Tudor},
	DOI = {10.1016/j.spa.2022.06.013},
	URL = {https://doi.org/10.1016/j.spa.2022.06.013},
}

@article {Balan_Yuan_HAMtimeindependentrough_2023,
	AUTHOR = {Balan, Raluca M. and Yuan, Wangjun},
	TITLE = {Hyperbolic {A}nderson model with time-independent rough noise:
	{G}aussian fluctuations},
	JOURNAL = {Electron. J. Probab.},
	FJOURNAL = {Electronic Journal of Probability},
	VOLUME = {29},
	YEAR = {2024},
	PAGES = {Paper No. 177, 52},
	ISSN = {1083-6489},
	MRCLASS = {60H15 (60F05 60G22 60H07)},
	MRNUMBER = {4836011},
	DOI = {10.1214/24-ejp1239},
	URL = {https://doi.org/10.1214/24-ejp1239},
}

@article {Kuzgun_Nualart_densitySHE_2022,
	AUTHOR = {Kuzgun, Sefika and Nualart, David},
	TITLE = {Convergence of densities of spatial averages of stochastic
	heat equation},
	JOURNAL = {Stochastic Process. Appl.},
	FJOURNAL = {Stochastic Processes and their Applications},
	VOLUME = {151},
	YEAR = {2022},
	PAGES = {68--100},
	ISSN = {0304-4149,1879-209X},
	MRCLASS = {60H15 (60H07)},
	MRNUMBER = {4441503},
	MRREVIEWER = {Xiaobin\ Sun},
	DOI = {10.1016/j.spa.2022.06.001},
	URL = {https://doi.org/10.1016/j.spa.2022.06.001},
}

@article {Kuzgun_Nualart_densitySHEcolored_2024,
	AUTHOR = {Kuzgun, Sefika and Nualart, David},
	TITLE = {Convergence of densities of spatial averages of the parabolic
	{A}nderson model driven by colored noise},
	JOURNAL = {Stochastics},
	FJOURNAL = {Stochastics. An International Journal of Probability and
	Stochastic Processes},
	VOLUME = {96},
	YEAR = {2024},
	NUMBER = {2},
	PAGES = {968--984},
	ISSN = {1744-2508,1744-2516},
	MRCLASS = {60H15 (60H07)},
	MRNUMBER = {4790990},
	DOI = {10.1080/17442508.2023.2238954},
	URL = {https://doi.org/10.1080/17442508.2023.2238954},
}

@article {Hu_Lu_Nualart_density_2014,
	AUTHOR = {Hu, Yaozhong and Lu, Fei and Nualart, David},
	TITLE = {Convergence of densities of some functionals of {G}aussian
	processes},
	JOURNAL = {J. Funct. Anal.},
	FJOURNAL = {Journal of Functional Analysis},
	VOLUME = {266},
	YEAR = {2014},
	NUMBER = {2},
	PAGES = {814--875},
	ISSN = {0022-1236,1096-0783},
	MRCLASS = {60F05 (60F10 60G15 60H07)},
	MRNUMBER = {3132731},
	MRREVIEWER = {Ivan\ Nourdin},
	DOI = {10.1016/j.jfa.2013.09.024},
	URL = {https://doi.org/10.1016/j.jfa.2013.09.024},
}

@article {Nualart_Zakai_1986,
	AUTHOR = {Nualart, David and Zakai, Moshe},
	TITLE = {Generalized stochastic integrals and the {M}alliavin calculus},
	JOURNAL = {Probab. Theory Relat. Fields},
	FJOURNAL = {Probability Theory and Related Fields},
	VOLUME = {73},
	YEAR = {1986},
	NUMBER = {2},
	PAGES = {255--280},
	ISSN = {0178-8051,1432-2064},
	MRCLASS = {60H07 (28C20 46G12 58G32 60G30 60H05)},
	MRNUMBER = {855226},
	MRREVIEWER = {Ana\ Bela\ Cruzeiro},
	DOI = {10.1007/BF00339940},
	URL = {https://doi.org/10.1007/BF00339940},
}

@book {Hu_BookGaussian_2017,
	AUTHOR = {Hu, Yaozhong},
	TITLE = {Analysis on {G}aussian Spaces},
	PUBLISHER = {World Scientific Publishing Co. Pte. Ltd., Hackensack, NJ},
	YEAR = {2017},
	PAGES = {xi+470},
	ISBN = {978-981-3142-17-6},
	MRCLASS = {60G15 (28C20 46-02 46E27 60B11 60H07 60Hxx)},
	MRNUMBER = {3585910},
	MRREVIEWER = {Jan\ van Neerven},
}

@article {Quer-Sardanyons_Sanz-Solé_3Dcontinuity_2004,
	AUTHOR = {Quer-Sardanyons, L. and Sanz-Sol\'e, M.},
	TITLE = {Absolute continuity of the law of the solution to the
	3-dimensional stochastic wave equation},
	JOURNAL = {J. Funct. Anal.},
	FJOURNAL = {Journal of Functional Analysis},
	VOLUME = {206},
	YEAR = {2004},
	NUMBER = {1},
	PAGES = {1--32},
	ISSN = {0022-1236,1096-0783},
	MRCLASS = {60H15},
	MRNUMBER = {2024344},
	MRREVIEWER = {Tom\'as\ Caraballo},
	DOI = {10.1016/S0022-1236(03)00065-X},
	URL = {https://doi.org/10.1016/S0022-1236(03)00065-X},
}
\end{document}